\documentclass[journal]{IEEEtran}
\usepackage{fancyhdr}
\usepackage[utf8]{inputenc}
\usepackage{verbatim}
\usepackage{amsmath}
\usepackage{amsthm}
\usepackage{amssymb}
\usepackage{graphicx}
\usepackage{lmodern}
\usepackage{mathdots}
\usepackage{graphicx}
\usepackage{amsmath}
\usepackage{amssymb}
\usepackage{graphicx}
\usepackage{mathpazo}
\usepackage{algorithm}
\usepackage{algorithmic}
\usepackage{array}
\usepackage{longtable}
\usepackage{multirow}

\usepackage[unicode=true,
 bookmarks=false,
 breaklinks=false,pdfborder={0 0 1},backref=section,colorlinks=false]
 {hyperref}
\hypersetup{
 colorlinks,bookmarksopen,bookmarksnumbered,citecolor=blue,urlcolor=red}

\makeatletter
 \let\oldforeign@language\foreign@language
 \DeclareRobustCommand{\foreign@language}[1]{%
   \lowercase{\oldforeign@language{#1}}}

 \let\oldforeign@language\foreign@language
 \DeclareRobustCommand{\foreign@language}[1]{%
   \lowercase{\oldforeign@language{#1}}}

\ifCLASSINFOpdf
\else
\fi

\newcommand{\Lone}{$ \mathcal{L}_1 $}

\newtheorem{rem}{Remark}

\newtheorem{assum}{Assumption}

\begin{document}
\fancyhf{}
\lhead{\small{To cite this article:
\bf{\textcolor{red}{Hashim A. Hashim, Sami El-Ferik, Babajide O. Ayinde, and Mohamed A. Abido. "Optimal Tuning of Fuzzy Feedback filter for L1 Adaptive Controller Using Multi-Objective Particle Swarm Optimization for Uncertain Nonlinear MIMO Systems." arXiv preprint arXiv:1710.05423 (2017).}}}}

\title{Optimal Tuning of Fuzzy Feedback filter for $ \mathcal{L}_1 $ Adaptive Controller Using Multi-Objective Particle Swarm Optimization for Uncertain Nonlinear MIMO Systems}

\author{Hashim~A.~Hashim, Sami~El-Ferik, Babajide O. Ayinde and~Mohamed~A.~Abido
\thanks{H. A. Hashim is with the Department of Electrical and Computer Engineering,
University of Western Ontario, London, ON, Canada, N6A-5B9, e-mail: \texttt{hmoham33@uwo.ca}  }
\thanks{S. El-Ferik is with the Department of Systems Engineering, King Fahd University of Petroleum and Minerals, Dhahran, 31261, Saudi Arabia e-mail: \texttt{selferik@kfupm.edu.sa}.}
\thanks{B. O. Ayinde is with the Department of Electrical and Computer Engineering, University of Louisville, Louisville, KY 40218, USA,  e-mail: \texttt{babajide.ayinde@louisville.edu}.}
\thanks{M. A. Abido is with Electrical and Computer Engineering Department, King Fahd University of Petroleum and Minerals, Dhahran, 31261, Saudi Arabia e-mail: \texttt{mabido@kfupm.edu.sa}.}}

\maketitle
\begin{abstract}
This paper proposes an efficient approach for tuning \Lone feedback filter of adaptive controller for multi-input multi-output (MIMO) systems. The feedback filter provides performance that trades off fast closed loop dynamics, robustness margin, and control signal range. Thus appropriate tuning of the filter's parameters is crucial to achieve optimal performance. For MIMO systems, the parameters tuning is challenging and requires a multi-objective performance indices to avoid instability. This paper proposes a fuzzy-based \Lone feedback filter design tuned with multi-objective particle swarm optimization (MOPSO) to remove these bottlenecks. MOPSO guarantees the appropriate selection of the fuzzy membership functions. The proposed approach is validated using twin rotor MIMO system and simulation results demonstrate the efficacy of here proposed while preserving the system stabilizability.
\end{abstract}



\IEEEpeerreviewmaketitle{}

\section{Introduction}

%
%
%
%
\IEEEPARstart{R}{ecently}, \Lone \,adaptive controller has been proposed to handle both single-input single-output (SISO) and multi-input multi-output (MIMO) nonlinear systems with uncertain parameters, unmodeled dynamics and/or unmeasurable external disturbances. \Lone \,adaptive controller provides fast adaptation and robustness to complete unknown dynamics. It was initially introduced for SISO system with unknown uncertainties \cite{cao_design_2006} and since then, it has been developed successfully for linear SISO systems with time varying uncertainties \cite{cao_design_2008}, nonlinear SISO systems with uncertainty \cite{cao_guaranteed_2007,luo_l_2010}, nonlinear MIMO systems with unmatched uncertainties \cite{xargay_l1_2010-1}, and for many other systems \cite{hovakimyan_l1_2010}. In addition to the aforementioned applications, \Lone \,adaptive controller has also been used extensively in aerospace applications \cite{wang_novel_2008,gregory_l1_2009}. It has three fundamental components namely the predictor, projection operators, and a low pass filter \cite{hovakimyan_l1_2010}. The low pass filter functions to mitigate the effect of both the system uncertainties and the frequency range of the control signal. Thus a careful design of this filter is crucial since it ensures both the fast adaptation and boundedness of transient and steady state performance.\\\\
\indent
The selection of the feedback filter coefficients has been long debated due to the trade-off between fast closed loop dynamics and robustness margin \cite{hovakimyan_l1_2010,kim_multi-criteria_2014,kim_filter_2011}. The optimal selection of the filter parameters for different structures have been studied extensively in \cite{hovakimyan_l1_2010}. Several heuristics have been proposed in the literature to estimate optimal filter coefficient. The list includes but not limited to using convex optimization with linear matrix inequalities (LMI) {\cite{li_filter_2008,hovakimyan_l1_2010}, MATLAB solver for multi-objective optimization \cite{li_optimization_2007}, and more recently a greedy randomized algorithm \cite{kim_multi-criteria_2014}.\\\\
\indent
It is remarked that the range and rate of the error between reference signals and actual outputs play an important role in the trade-off between fast closed loop dynamics and robustness margin. Most of the previous work that address filter design assumed constant filter coefficients \cite{hovakimyan_l1_2010,kim_multi-criteria_2014,kim_filter_2011,li_filter_2008,li_optimization_2007,kharisov_limiting_2011} and optimal filter parameters are selected offline. One important inference from these work highlights the tradeoff between filter bandwidth and the robustness margin and to avoid deteriorating the robustness, slower closed loop dynamics is imposed. However, this imposition degrades the output performance of the system \cite{hovakimyan_l1_2010}. This problem is more pronounced especially in the case of trajectory tracking problems. One way to overcome this limitation is
to tune feedback filter parameters online in order to into account the adverse relationships between error values and the feedback filter gains \cite{hashim_fuzzy_2015}.\\\\
\indent
Fuzzy logic controller (FLC) has been widely used in control applications to design smooth control signals for nonlinear systems. It has also been used to tune controller parameters for enhancing closed loop performance. For instance, FLC has been deployed for optimal tuning of PID controllers \cite{pan_tuning_2011,sharma_performance_2014}. FLC achieves optimality in many applications through proper choice of membership functions, scaling factors of input values, and rule-based \cite{reznik_pid_2000}. Defining membership functions, linguistic variables, and other parameters of FLC using a trial and error approach is nontrivial and time consuming. To circumvent this bottleneck, global optimization tools such as Particle Swarm Optimization (PSO) can be employed to find optimal solutions in the search space with predefined constraints.\\\\
\indent
PSO \cite{eberhart_new_1995} is a global search technique that has been found efficient in many control applications and other fields of study. It has been used to tune membership functions variables \cite{hashim_fuzzy_2015}, and find optimal control parameters of adaptive fuzzy controller \cite{das_sharma_hybrid_2009}. In addition, the multi objectives version of PSO known as Multi-Objective Particle Swarm Optimization (MOPSO) was introduced to obtain best compromise solution among many conflicting objectives \cite{abido_multiobjective_2009}. In total, heuristic techniques provide effective  global search solution for complex problems \cite{hashim_fuzzy_2015, abido_multiobjective_2009, Abd_2017}.\\\\
\indent
This paper considerably expands the scope of the fuzzy-\Lone \,adaptive controller for nonlinear SISO system introduced in \cite{mohamed_improved_2014,hashim_fuzzy_2015}. In \cite{hashim_fuzzy_2015}, PSO is responsible for finding best parameters of the input-output membership function based on a trade-off cost functions:\\
   \begin{equation}
                  \label{eq:chFuzL1Obj1}
                  \begin{aligned}
                   Obj = \sum\limits_{t=0}^{t_{\rm sim}}\big(\gamma_1 e^2\left(t\right)+\gamma_2 u^2\left(t\right)\big)
                  \end{aligned}
                \end{equation}
where $e\left(t\right) = r\left(t\right) - y\left(t\right)$, $e\left(t\right)$ and $u\left(t\right)$ are the system error and control signal respectively. $\gamma_1$ and $\gamma_2$ are weights that can be selected arbitrarily. The single objective function in \eqref{eq:chFuzL1Obj1} is designed to handle SISO systems which could result in loss of controllability when directly applied to MIMO systems. A closer scrutiny of \eqref{eq:chFuzL1Obj1} reveals that one objective could be favored over the other because the implementation results in linear weighting of two objectives. By consequence, the imposed weights and the final solution confines the solution on the pareto front and will not guarantee a best compromise solution. To avoid this problem, a Multi-objective PSO (MOPSO) is here proposed to optimize the parameters of input and output membership functions of FLC-based $\mathcal{L}_1$ \,adaptive controller. Fuzzy logic is implemented to tune the feedback filter gains according to the values of the error and its rate. A compromise solution can be obtained between control signal and error regardless of the difference in their range. The best compromise solution will be obtained through the Pareto-optimal front which is regarded to be a more feasible solution \cite{abido_multiobjective_2009}. It is remarked that during preliminary experiment for small initial error and in accordance with robustness margin reduction, traditional $\mathcal{L}_1$ \,adaptive controller became unstable. This observation reinforces one of the limitations of traditional $\mathcal{L}_1$ \,adaptive controller. For this reason, we could not benchmark the proposed enhancement with traditional $\mathcal{L}_1$ \,adaptive controller in the result discussion section. In the future, the problem can be extended to regulate the output consensus of heterogeneous uncertain nonlinear multi-agent systems for example (\cite{Hashim2017adaptive,el2017neuro,hashim2017neuro}), using tuned $\mathcal{L}_1$ \,adaptive controller.\\\\
\indent
The problem addressed here is two-fold: $\left(i\right)$ fuzzy-\Lone \,adaptive controller is developed for nonlinear MIMO systems, (ii) multi-objective optimization is adapted to find the optimal parameters of the tuned filter. This paper is organized as follows. Section \ref{Sect2} presents a brief overview of \Lone \,adaptive controller for unmatched nonlinear systems. Section \ref{Sec3} details the idea of fuzzy filter design and the structure of the proposed control. Section \ref{Sec4} formulates the optimization problem and presents multi-objective particle swarm optimization for FLC design. Illustrative examples are used to validate the robustness of the proposed approach in section \ref{Sec5}. Concluding remarks are highlighted in section \ref{Sec6}.
\section{General Overview of \Lone \,adaptive controller}\label{Sect2}
For convenience, this section gives a brief review of \Lone \,adaptive control design. Consider the following nonlinear system dynamics
    \begin{equation}
     \begin{aligned}
      \label{eq:L1Act1}
        \dot{x}\left(t\right) = & A_{m}x\left(t\right) + B_{m}\omega u\left(t\right) + f\left(x\left(t\right),z\left(t\right),t\right),\hspace{2pt}x\left(0\right)=x_0\\
        \dot{x}_z\left(t\right) = & g(x\left(t\right),x_{z}\left(t\right),t),\hspace{10pt}x_{z}\left(0\right)=x_{z_0}\\
        z\left(t\right) = & g_{0}\left(x_{z}\left(t\right),t\right)\\
        y\left(t\right) = & Cx\left(t\right)
     \end{aligned}
    \end{equation}
where $x\left(t\right) \in \mathbb{R}^{n}$ is the system measured state vector, $u\left(t\right) \in \mathbb{R}^{m}$ is the control input vector, $y\left(t\right) \in \mathbb{R}^{m}$ is the system output vector, $B_{m} \in \mathbb{R}^{n \times m}$ is the desired input state matrix which is assumed to be known and constant with full column rank; the pair $\left(A_m,B_{m}\right)$ is controllable and $C \in \mathbb{R}^{m \times n}$ is known output state matrix and constant with full row rank and $\left(A_m,C_{m}\right)$ is observable. $A_m \in \mathbb{R}^{n \times n}$, known as Hurwitz matrix, includes the desired dynamics for the closed-loop system, $\omega \in \mathbb{R}^{m \times m}$ is a gain matrix that signifies uncertain system input; $x_{z}\left(t\right)$ is the state vector of unmodeled internal dynamics and $z\left(t\right)$ is the output of internal dynamics. $f:\mathbb{R}\times\mathbb{R}^{n}\times\mathbb{R}^{p}\to\mathbb{R}^{n}$, $g_0:\mathbb{R}^{l}\times\mathbb{R}\to\mathbb{R}^{p}$ and $g:\mathbb{R}\times\mathbb{R}^{l}\times\mathbb{R}^{n}\to\mathbb{R}^{l}$ are unknown nonlinear continuous functions.\\\\
The system in (\ref{eq:L1Act1}) can be written as
    \begin{equation}
     \begin{aligned}
      \label{eq:L1Act2}
       \dot{x}\left(t\right) =  & A_{m}x\left(t\right) + B_{m}\left(\omega u\left(t\right) + f_{1}\left(x\left(t\right),z\left(t\right),t\right)\right) \\
       & + B_{um}\left(f_2\left(x\left(t\right),z\left(t\right),t\right)\right) ,\hspace{10pt}x\left(0\right)=x_0\\
       \dot{x}_z =  & g\left(x\left(t\right),x_{z}\left(t\right),t\right),\hspace{10pt}x_{z}\left(0\right)=x_{z_0}\\
        z\left(t\right) =  & g_{0}\left(x_{z}\left(t\right),t\right)\\
        y\left(t\right) =  & Cx\left(t\right)
     \end{aligned}
    \end{equation}
With reference to (\ref{eq:L1Act2}), highly nonlinear system with strong coupling and unmatched nonlinearities is divided into two parts. The first part $f_{1}\left(\cdot\right)$ includes the matched unknown nonlinear components and the second part $f_2\left(\cdot\right)$ contains unmatched unknown nonlinear components. $B_{um} \in \mathbb{R}^{n \times \left(n-m\right)}$ is a constant matrix and should be selected such that $B_{m} \times B_{um}=0$ and $rank([B_{m},B_{um}])=n$. Let $X \triangleq [x^{\top},z^{\top}]$, $f_{1}\left(t,X\right):\mathbb{R}\times\mathbb{R}^{n}\times\mathbb{R}^{p}\to\mathbb{R}^{m}$ and $f_2\left(t,X\right):\mathbb{R}\times\mathbb{R}^{n}\times\mathbb{R}^{p}\to\mathbb{R}^{n-m}$ be unknown nonlinear functions that align with the following assumptions \cite{xargay_l1_2010-1,hovakimyan_l1_2010}.

\begin{assum} The control input is partially known with known sign and the system input gain matrix $\omega$ is assumed to be nonsingular and unknown with strictly row-diagonally dominant matrix form but with the signs of diagonal elements known.
\begin{equation*}
        \omega \in \Omega \subset \mathbb{R}^{m \times m}
\end{equation*}
where $\Omega$ is assumed to be known convex compact set.
\end{assum}

\begin{assum} Let $B\in \mathbb{R}^+$, $f\left(0,t\right)$ be uniformly bounded such that $ f_i\left(0,t\right) \leq B $ $\forall$ $ t \geq 0 $
\end{assum}

\begin{assum} Partial derivatives of the nonlinear functions are continuous and uniformly bounded, where for any $\delta >0$, there exist $d_{f_{xi}} \left(\delta\right) >0 $ and $d_{f_{ti}} \left(\delta\right) >0 $ such that for arbitrary $\left\Vert x \right\Vert_{\infty} \leq \delta$ and any $u$, the partial derivatives of $ f_i\left(t,X\right)) $ is piecewise-continuous and bounded,
\begin{equation*}
       \left\Vert \frac{\partial f_i\left(t,X\right) }{\partial x} \right\Vert \leq d_{f_{xi}} \left(\delta\right),\hspace{10pt} \left\Vert \frac{\partial f_i\left(t,X\right) }{\partial t} \right\Vert \leq d_{f_{ti}} \left(\delta\right) \hspace{10pt}i=1,2
\end{equation*}
\end{assum}

\begin{assum} The internal dynamics are BIBO stable with respect to $x_{z0}$ and $x\left(t\right)$ and there exist $L_z>0$ and $B_z>0$ such that
\begin{equation*}
        \left\Vert z_t \right\Vert_{\mathcal{L}_{\infty}} \leq L_z\left\Vert x\left(t\right) \right\Vert_{\mathcal{L}_{\infty}}+B_z \forall t \geq 0
\end{equation*}
\end{assum}

\begin{assum} Transmission zeros of the transfer matrix are stable where zeros of $H_m\left(s\right) = C(sI-A_m)^{-1}B_{m}$ are located in the open left half of the complex plane.
\end{assum}

Considering the following notations:
\begin{equation*}
       H_{xm}\left(s\right) \triangleq \left(s\mathbb{I}_{n}-A_{m}\right)^{-1}B_{m},
\end{equation*}
     \begin{equation*}
       H_{xum}\left(s\right) \triangleq \left(s\mathbb{I}_{n}-A_{m}\right)^{-1}B_{um},
     \end{equation*}
     \begin{equation*}
       H_{m}\left(s\right) \triangleq C\left(s\mathbb{I}_{n}-A_{m}\right)^{-1}B_{m},
     \end{equation*}
     \begin{equation*}
       H_{um}\left(s\right) \triangleq C\left(s\mathbb{I}_{n}-A_{m}\right)^{-1}B_{um},
     \end{equation*}
$D\left(s\right)$ is a strictly proper transfer matrix,  $K \in \mathbb{R}^{m \times m}$ with $\omega \in \Omega$ included in \Lone \,adaptive controller and guarantees stability of the strictly proper transfer function $C\left(s\right)$ such as:
     \begin{equation}
       C\left(s\right) \triangleq \omega K D\left(s\right)\left(\mathbb{I}_{m}+\omega K D\left(s\right)\right)^{-1}
     \end{equation}
     $D\left(s\right)$ should be selected such that $C\left(s\right)H^{-1}\left(s\right)$ is a proper stable transfer matrix and in this case, $D\left(s\right)=1/s \cdot \mathbb{I}_m$, which leads to
     \begin{equation}
       C\left(s\right) \triangleq \omega K\left(s\mathbb{I}_{m}+\omega K\right)^{-1}
     \end{equation}
The structure of \Lone \,adaptive controller can simply be described as in Fig. \ref{Fuzzy_L1_Gen} where the controller is divided into state predictor, projection operators and feedback filter. Although \Lone \,adaptive controller allows decoupling between adaptation and robustness margin through high gain, the structure of \Lone \,adaptive controller introduces coupling between fast closed loop dynamics and robustness margin.
     \newline{\bf State Predictor:} The following state predictor is considered:
     \begin{equation}
       \label{eq:L1est}
       \begin{aligned}
          &\begin{split}
          \dot{\hat{x}}\left(t\right) = & A_{m}\hat{x}\left(t\right)  +  B_{m}\left(\hat{\omega} u\left(t\right) + \hat{\theta}_1 \left\Vert x\left(t\right) \right\Vert_{\infty}+\hat{\sigma}_1 \right)\\
          &+ B_{um}\left(\hat{\theta}_2 \left\Vert x\left(t\right) \right\Vert_{\infty}+\hat{\sigma}_2 \right),\hspace{10pt}\hat{x}\left(0\right) = x\left(0\right)
       \end{split}\\
         & \hat{y}\left(t\right) = c\hat{x}\left(t\right)
       \end{aligned}
     \end{equation}
where $\hat{x} \in \mathbb{R}^{n}$ is the predicted state vector and $\hat{y} \in \mathbb{R}^{m}$ is the predicted output vector. $\hat{\omega} \in \mathbb{R}^{m \times m}$, $\hat{\theta}_1\left(t\right)\in \mathbb{R}^{m}$, $\hat{\theta}_2\left(t\right)\in \mathbb{R}^{n-m}$, $\hat{\sigma}_1\left(t\right)\in\mathbb{R}^{m}$ and $\hat{\sigma}_2\left(t\right)\in\mathbb{R}^{n-m}$ are all adaptive estimates and they are defined using the following adaptation laws \cite{hovakimyan_l1_2010}.
     \begin{equation}
       \label{eq:L1Proj}
       \begin{aligned}
         &\dot{\hat{\omega}} = \Gamma {\rm Proj}\left(\hat{\omega},-\left(\tilde{x}^{\top}PB_{m}\right)^{\top}\left(u\left(t\right)^{\top}\right)\right), \hat{\omega}\left(0\right) = \hat{\omega}_0  \\
         &\dot{\hat{\theta}}_1 = \Gamma {\rm Proj}\left(\hat{\theta}_1,-\left(\tilde{x}^{\top}PB_{m}\right)^{\top}\left\Vert x\left(t\right) \right\Vert_{\infty}\right), \hat{\theta}_1\left(0\right) = \hat{\theta}_{10}\\
         &\dot{\hat{\theta_2}} = \Gamma {\rm Proj}\left(\hat{\theta}_2,-\left(\tilde{x}^{\top}PB_{um}\right)^{\top}\left\Vert x\left(t\right) \right\Vert_{\infty}\right), \hat{\theta}_2\left(0\right) = \hat{\theta}_{20}\\
         &\dot{\hat{\sigma}}_1 = \Gamma {\rm Proj}\left(\hat{\sigma}_1,-\left(\tilde{x}^{\top}PB_{m}\right)^{\top}\right), \hat{\sigma}_1\left(0\right) = \hat{\sigma}_{10}\\
         &\dot{\hat{\sigma}}_2 = \Gamma {\rm Proj}\left(\hat{\sigma}_1,-\left(\tilde{x}^{\top}PB_{m}\right)^{\top}\right), \hat{\sigma}_2\left(0\right) = \hat{\sigma}_{20}
       \end{aligned}
     \end{equation}
where $\tilde{x} \triangleq \hat{x} - x\left(t\right)$, $\Gamma \in \mathbb{R}^{+}$ is the adaptation gain, and $P$ is the solution of Lyapunov equation $A_m^{\top}P+PA_m=-Q$ with $P$ and $Q$ are symmetric and positive definite matrices. The projection operator ensures that $\hat{\omega} \in \Omega$, $||\hat{\theta}_i||_{\infty} \in \Theta_i$, $||\hat{\sigma}_i|| \leq \Delta_i$, where $\theta_{bi}$ and $\delta_{bi}$ are determined numerically. Projection operators will be evaluated as defined in \cite{pomet_adaptive_1992}.\\\\
{\bf Control Law:} Control signal is calculated as follows
     \begin{equation}
       u\left(s\right) = -kD\left(s\right)\hat{\eta}\left(s\right)
     \end{equation}
where $r\left(s\right)$ and $\hat{\eta}\left(s\right)$ are the Laplace transforms of $r\left(t\right)$ and $\hat{\eta}\left(t\right) = \hat{\omega} u\left(t\right) + \hat{\eta}_1 + \hat{\eta}_2 - K_gr\left(t\right)$ respectively. The feedforward gain for zero steady state error is calculated using $K_g \triangleq -1(CA_{m}^{-1}B_{m})^{-1}$; $K > 0$ is a feedback diagonal matrix gain and $D\left(s\right) = \frac{1}{s}$ is a strictly proper transfer function. $D\left(s\right)$ and $K$ ensure strictly proper stable closed loop system with $\hat{\eta}_1$ and $\hat{\eta}_2$ evaluated using
     \begin{equation}
       \hat{\eta}_1 \triangleq \hat{\theta}_1 \left\Vert x\left(t\right) \right\Vert_{\infty}+\hat{\sigma}_1
     \end{equation}
     \begin{equation}
       \hat{\eta}_2 \triangleq \hat{\theta}_2 \left\Vert x\left(t\right) \right\Vert_{\infty}+\hat{\sigma}_2
     \end{equation}
The DC gain is $C\left(0\right) = \mathbb{I}_m$. More details on \Lone \,adaptive controller for highly nonlinear unmatched system are given in \cite{xargay_l1_2010-1,hovakimyan_l1_2010}. The schema of \Lone \,adaptive control is depicted in Fig. \ref{Fuzzy_L1_Gen}.

    \begin{figure*}[ht]
             \centering
             \includegraphics[scale=0.8]{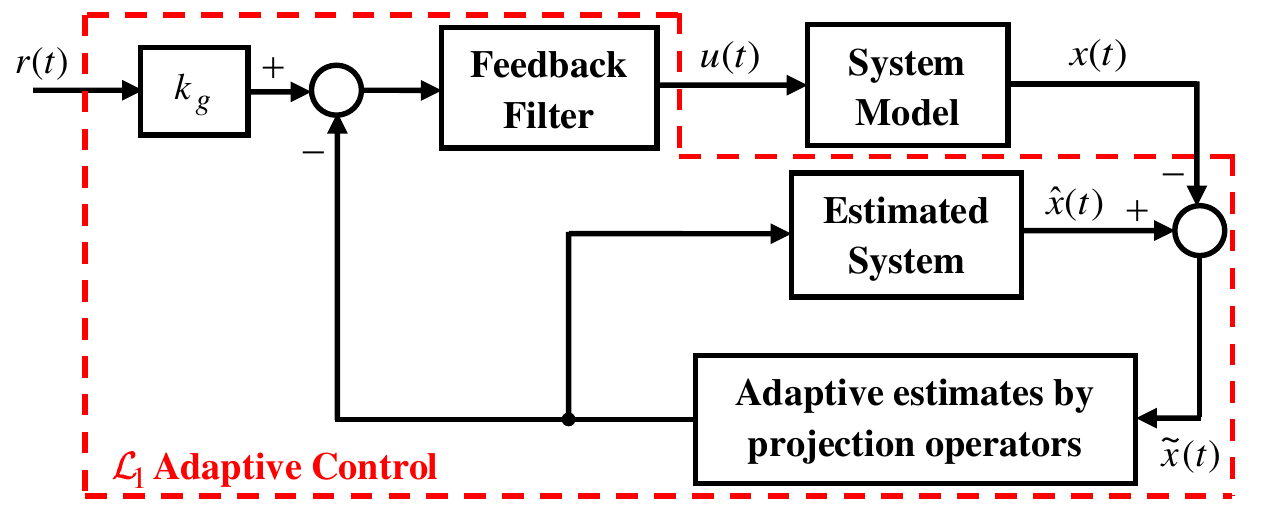}
             \caption{General structure of \Lone \,adaptive controller }\label{Fuzzy_L1_Gen}
     \end{figure*}

\section{Optimal Fuzzy-tuning of the feedback filter}\label{Sec3}
In this work, the main objectives are to design FLC to tune the parameters of \Lone \,adaptive feedback filter and to produce a smooth output signal $y\left(t\right)$ that tracks a reference signal $r\left(t\right)$ with the desired performance as depicted in Fig.~\ref{Fuzzy_L1_3PS}. The tuned filter enables the selection of fast closed loop dynamics with proper robustness margin. The appropriate parameters of the feedback filter are determined online during the control process.
\subsection{Structure of Fuzzy Logic Controller}
The difference between the regulated output vector $y\left(t\right)$ and the reference input vector $r\left(t\right)$ is the error vector $e\left(t\right)$. The infinity norm of error $e\left(t\right)$ and rate of error $\dot{e}\left(t\right)$ are the two fuzzy inputs. Each of $\left\Vert e\left(t\right) \right\Vert_{\infty}$ and $||\dot{e}\left(t\right)||_{\infty}$ are multiplied by weight gains $k_p$ and $k_d$ respectively, where $k_p$ and $k_d$ are proportional and the differential weights respectively. The selection of these gains will be adjusted such that the input of FLC is normalized between 0 and 1.
\begin{equation}
    k_p \leq \frac{1}{\left\Vert e \right\Vert_{\infty} }, \hspace{10pt} k_d \leq \frac{1}{\left\Vert \dot{e} \right\Vert_{\infty}}
\end{equation}
The selection of the norms guarantees stable dynamics of \Lone \,adaptive controller with FLC feedback filter. The output of FLC  is the inverse of the feedback filter gain $K_{i,i} = 1/k_f$ where $i=1,\ldots , m$. The feedback gain matrix of \Lone \,adaptive controller will be selected such as $K = k\mathbb{I}_m$ if $\left\Vert e\left(t\right) \right\Vert_{\infty}$ is less than or equal $k_e$ where $k$ is a constant.  $k_e$ is a constant value and will be defined based on prior knowledge. The formulation of MOPSO will be covered in details in Section 3. Fig. \ref{Fuzzy_L1_3PS} gives the general schema of the proposed control structure.
    \begin{figure*}[ht]
             \centering
             \includegraphics[scale=0.8]{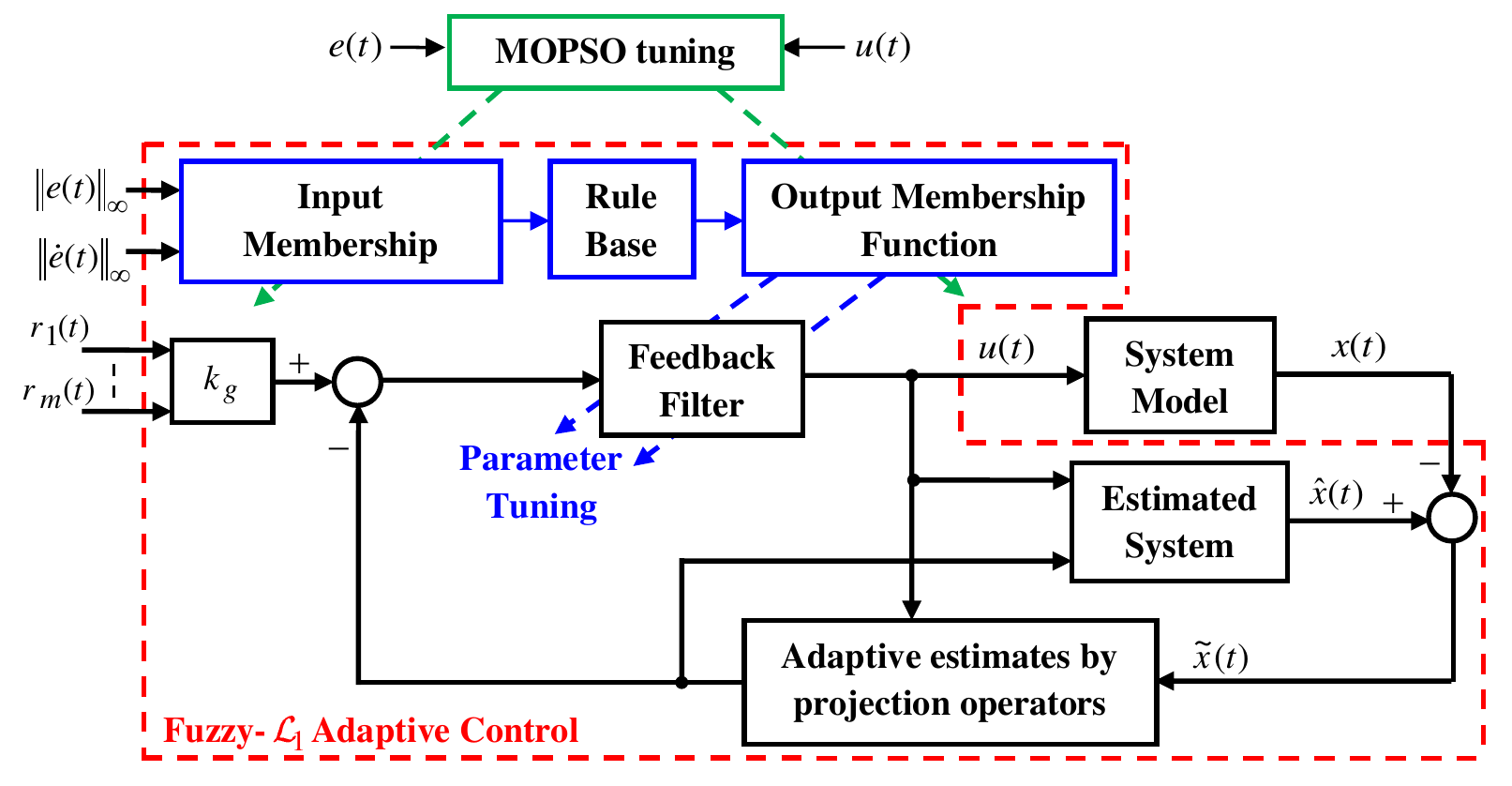}
             \caption{Proposed fuzzy-\Lone \,adaptive control structure with MOPSO Tuning.}
             \label{Fuzzy_L1_3PS}
          \end{figure*}

In \Lone \,adaptive controller, fast closed loop dynamics improve tracking capabilities but increase the control signal range and reduce robustness margin. Fuzzy-\Lone \,adaptive controller is proposed to ensure fast closed loop performance and to enhance the robustness margin. The design of FLC considers two objective functions that account for the control signal range and summation of the tracking error. Reducing the control signal range contradicts the reduction of error tracking. Therefore, multi-objective optimization technique is necessary to achieve an optimal compromise solution. In this scenario, a set of trade-off solutions will be obtained. This set is otherwise known as Pareto-optimal front \cite{zitzler_evolutionary_1998}. The input and output of FLC membership functions will be optimized using MOPSO.

\section{Multi-Objective Particle Swarm Optimization}\label{Sec4}
Particle Swarm Optimization (PSO) is an evolutionary heuristic that mimics the social behavior of bird swarming or fish schooling \cite{eberhart_new_1995}. PSO initiates population of particles randomly in space with each particle representing a potential solution. Each particle has a set of parameters and moves randomly in a multi-dimensional space in search of optimal solution. The velocity of each particle in space has a significant role in targeting the best candidate solution. In addition, velocity and position adjustments for each particle rely on the experiences gained from its own velocity, location and neighboring particles' locations. The velocity and position of each particle are updated as follows:
          \begin{equation}
          \label{eq:FuzL1v}
          \begin{split}
           v_{i,d}\left(t\right) = & \alpha\left(t\right)v_{i,d}\left(t-1\right) + c_{1}r_{1}\left(p_{i,d}^{*}\left(t-1\right)-p_{i,d}\left(t-1\right)\right) \\
           &+ c_{2}r_{2}\left(p_{i,d}^{**}\left(t-1\right)-p_{i,d}\left(t-1\right)\right)
          \end{split}
          \end{equation}
          \begin{equation}
          \label{eq:FuzL1x}
           p_{i,d}\left(t\right) = v_{i,d}\left(t\right) + p_{i,d}\left(t-1\right)
          \end{equation}
where $P_i$ is the candidate solution with $P_i\left(t\right) = [p_{i,1}\left(t\right),\ldots,p_{i,M}\left(t\right)] \in \mathbb{R}^M$; $M$ is the number of optimized parameters, $V_i$ is the velocity of candidates given as $V_i\left(t\right) = [v_{i,1}\left(t\right),\ldots,v_{i,M}\left(t\right)] \in \mathbb{R}^M$, $i=1,2,\ldots, N$; $N$ is the population size, and $P\left(t\right) = [P_{1}\left(t\right),\ldots,P_{N}\left(t\right)] \in \mathbb{R}^{N\times M}$. It should also be noted that across all dimensions the velocity should be bounded such that $v_{i,d}\left(t\right)\in [-v_{i,d}^{\rm max}, v_{i,d}^{\rm max}]$, where the maximum velocity is defined as in \eqref{eq:FuzL1x1} \cite{abido_multiobjective_2009}.
          \begin{equation}
          \label{eq:FuzL1x1}
           v_{i,d}^{\rm max} = \frac{x_{i,d}^{\rm max}-x_{i,d}^{\rm min}}{N_{int}}
          \end{equation}
where $N_{int}$ is the number of intervals and $d=1,2,\ldots, M$; $p_{i}^{*}$ and $p_{i}^{**}$ are local and global best solutions for each particle respectively. Other parameter settings such as $c_1$ and $c_2$ are personal and social behavior of parameters and $r_1$ and $r_2$ are randomly set to values between 0 and 1 \cite{eberhart_new_1995}.\\

In MOPSO, non-dominated local best set $S_{i}^{*}$ with pre-specified size stores a set of non-dominated solutions. At the initial stage, non-dominated local set starts with $S_{i}^{*} = P_i\left(0\right)$ after which non-dominated solutions will be added to the set with predefined size. Clustering is employed to reduce size of the non-dominated local set to a predefined value. Average distance between two pairs of clusters will be evaluated and the minimal distance will be combined into one cluster. Larger distances are retained for search enhancement and for coverage of more space. Non-dominated global set $S_{i}^{**}$ stores all non-dominated solutions starting from $P_i\left(0\right)$ up to last iteration. Similarly, clustering algorithm is implemented to reduce the non-dominated global set into a predefined set size. The output of non-dominated global set clustering will constitute the Pareto optimal front. All historical records of non-dominated solutions through the search process is stored in the external. The external set is updated continuously through the dominance algorithm \cite{abido_multiobjective_2009} and then the clustering algorithm is subsequently used to find the non-dominated solutions of the union between external and global set.\\
Local best $P_{i}^{*}$ and global best $P_{i}^{**}$ respectively belong to $S_{i}^{*}$ and $S_{i}^{**}$. The complete multi-objective algorithm can be found in \cite{abido_multiobjective_2009} and the multi-objective optimization problem can be formulated as follows
             \begin{equation}
               \label{eq:chFuzL1Obj}
               \begin{aligned}
                & E\left(i\right) = \sum\limits_{t=0}^{t_{\rm sim}}\big[ e_{1}^2\left(t\right)+\ldots+e_{m}^2\left(t\right)\big]\\
                & U\left(i\right) = \left\Vert u_{1}\left(t\right) \right\Vert_{\infty}+\ldots+\left\Vert u_{m}\left(t\right) \right\Vert_{\infty}\\
                & \min Obj(i,:) = \big(E\left(i\right), U\left(i\right) \big)\\
                & p_{i,d}^{\rm min}\ \leq p_{i,d}\left(t\right) \leq p_{i,d}^{\rm max}\\
                & v_{i,d}^{\rm min}\ \leq v_{i,d}\left(t\right) \leq v_{i,d}^{\rm max}\\
               \end{aligned}
             \end{equation}
where $P_i$ is the candidate solution that contains the optimized parameters or variables of the membership functions; $e\left(t\right) = r\left(t\right) - y\left(t\right)$; $e\left(t\right)$ and $u\left(t\right)$ are system's error and control signal respectively. Parameters of input and output membership functions are optimized to minimize the constrained objective function in (\ref{eq:chFuzL1Obj}). MOPSO is formulated to search for the Pareto-optimal front and a compromise solution will be selected from the set of nondominated solutions.\\
The flow diagram of MOPSO algorithm is shown in Fig. \ref{Fuzzy_L1_PSO}. The algorithm is used with \Lone \,adaptive controller for uncertain MIMO nonlinear systems to find the optimal parameters of fuzzy membership functions for a predefined number of generations \cite{abido_multiobjective_2009}.
       \begin{figure}[!h]
          \centering
          \includegraphics[scale=0.6]{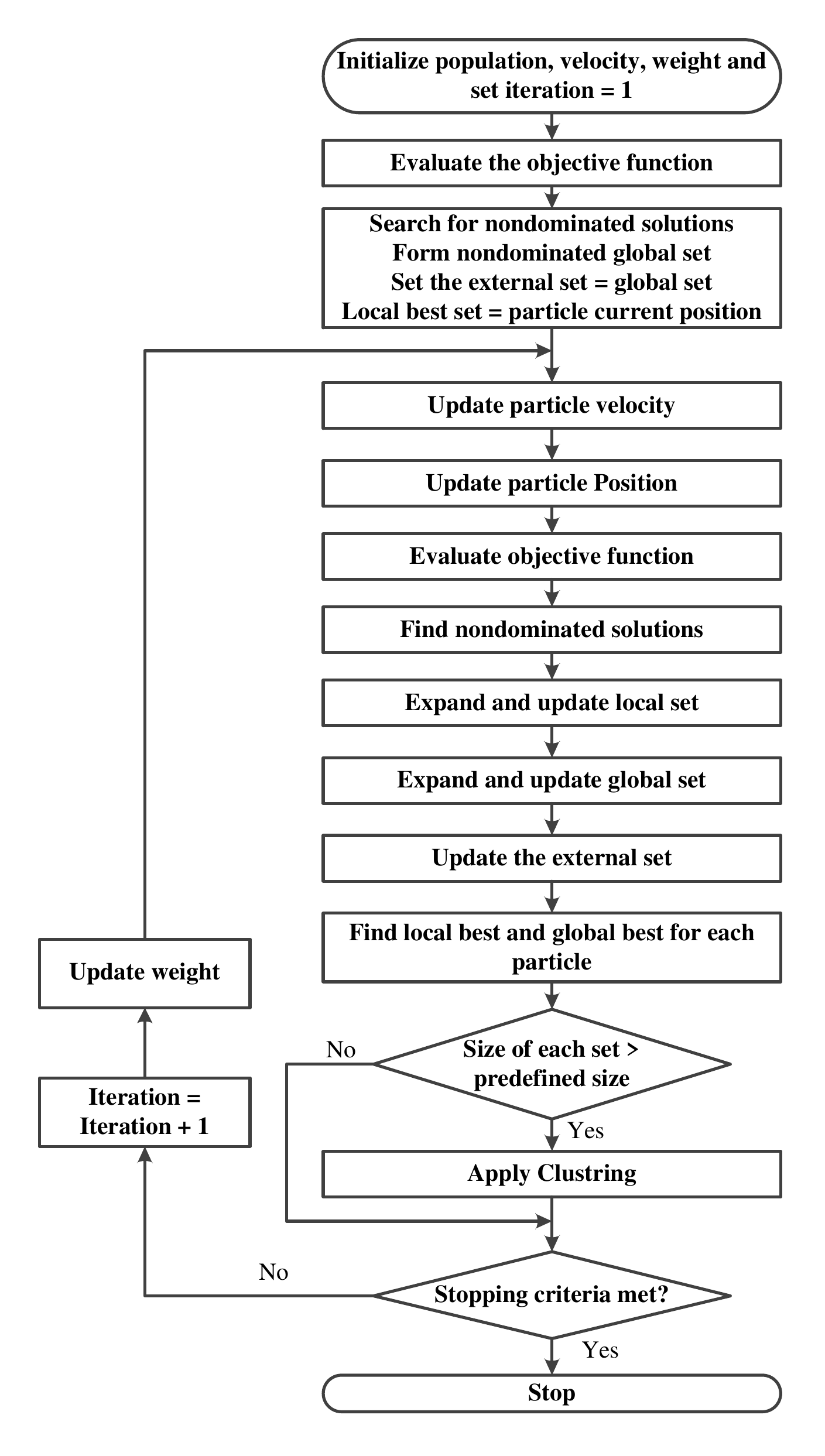}
          \caption{Flowchart of multi-objective particle swarm Optimization.}\label{Fuzzy_L1_PSO}
       \end{figure}

\begin{rem}
In the proposed approach filter properties such as strictly proper and low pass characteristics with $C\left(0\right)=1$ are preserved. In addition, the control input constitutes an independent objective in the optimization that freely moves within the sets. Controllability and stability of the Fuzzy-based-\Lone \,adaptive controller are preserved in line with stability analysis in \cite{cao_design_2006}.
\end{rem}

\section{Results and Discussions}\label{Sec5}

The performance in terms of tracking capability and robustness of the proposed controller is evaluated using a highly nonlinear unmatched system with strongly coupled dynamics. Robustness of the controller was examined by imposing uncertainties in the TRMS parameters. TRMS was chosen to evaluate the performance of \Lone\,adaptive control because it belongs to the class of unmatched system with very aggressive model nonlinearity and coupled dynamics. Also, the system is nonlinear in terms of control input $u$.
\subsection{Twin Rotor MIMO System Description}
Twin Rotor has strong coupling between the main and tail rotors and it emulates the helicopter dynamics in the pitch and yaw angle dynamics \cite{TWIN_twin_2002,Hash:2015}. The laboratory set up of TRMS is shown in Fig. \ref{TRMS_SET}.\\
       \begin{figure}[!h]
          \centering
          \includegraphics[scale=0.25]{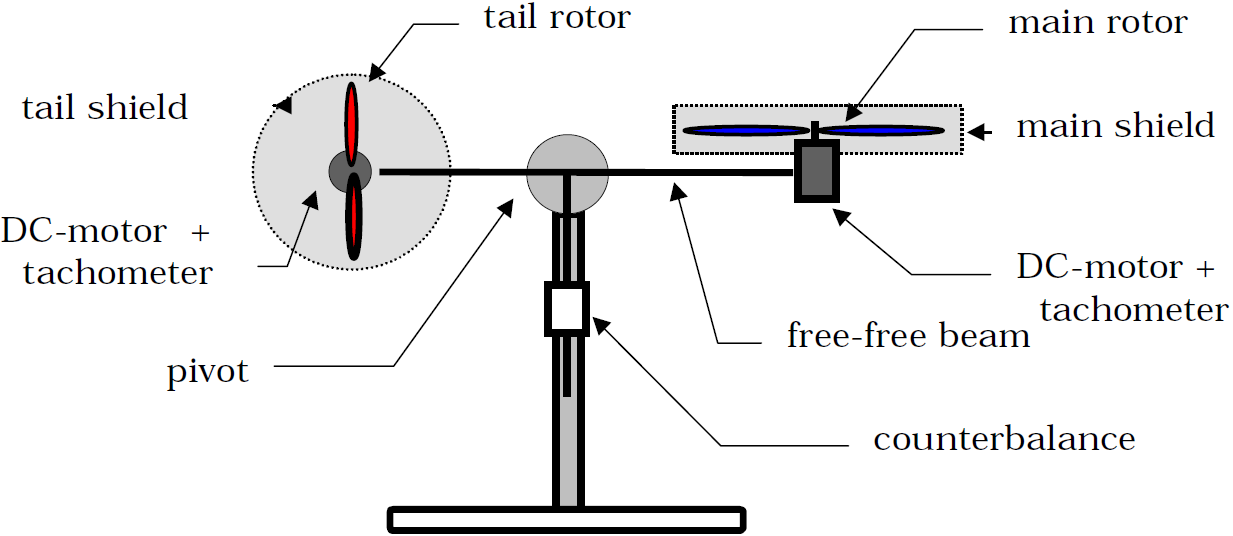}
          \caption{Laboratory set-up of TRMS.}\label{TRMS_SET}
       \end{figure}
The system is controlled by two control signals $u_1$, $u_2$ which are transferred into momentum torques $\tau_{1}$ and $\tau_2$. The pitch angle $\psi$ describes the motion of the main propeller which is to be controlled vertically. The tail rotor is controlled in the horizontal direction through the yaw angle $\phi$. Pivot beam with a weighted mass in TRMS is installed for stabilization and to allow free motions.\\
The mathematical model of TRMS is described by six states: vertical angle, yaw angle, pitch angular velocity, yaw angular velocity, and two momentum torques. Two potentiometers are fitted at the pivot in order to measure $\phi$ and $\psi$. The full description of TRMS are detailed in \cite{TWIN_twin_2002}. The model of TRMS is as given in \cite{mohamed_improved_2014,Hash:2015}:
     \begin{equation*}
      \begin{aligned}
       \label{eq:model}
        \dot{x_1} =& x_2\\
        \dot{x_2} =& \frac{a_3}{a_1}x_5^2 + \frac{a_5}{a_1}x_5 -\frac{a_7}{a_1}{\rm sin}(x_1) - \frac{a_8}{a_1}x_2 + \frac{0.0362}{a_1}x_4^2{\rm sin}(2x_1)\\
        & - \frac{a_{12}}{a_1}\big(a_3x_5^2 + a_5x_5\big)x_4{\rm cos}(x_1)\\
        \dot{x_3} =& x_4\\
        \dot{x_4} =& \frac{a_4}{a_2}x_6^2 + \frac{a_6}{a_2}x_6 -\frac{a_{10}}{a_2}x_4 - 1.75(k_c/a_2)\big(a_3x_5^2 + a_5x_5\big)\\
        \dot{x_5} =& -\frac{a_{16}}{a_{15}}x_5 + \frac{a_{13}}{a_{15}}u_1\\
        \dot{x_6} =& -\frac{a_{18}}{a_{17}}x_5 + \frac{a_{14}}{a_{17}}u_2\\
      \end{aligned}
     \end{equation*}
where the state vector $[\psi,\dot{\psi},\phi,\dot{\phi},\tau_{1},\tau_2]^{\top}$ are designated as $x=[x_1,x_2,x_3,x_4,x_5,x_6]^{\top}$; the output vector $y = [\psi,\phi]^{\top}$ as $[x_1,x_3]^{\top}$ and input vector is $u=[u_1,u_2]^{\top}$. The parameters of TRMS used for the experiments are listed in   \nameref{App_Sec6}, Table \ref{table:Tab_Twin}.
\subsection{Fuzzy-\Lone \,adaptive controller implementation:}
The TRMS model in more general form can be expressed as:
     \begin{equation*}
      \begin{aligned}
        & \begin{split}
        \dot{x}\left(t\right) = & A_{m}x\left(t\right) + B_{m}(\omega u\left(t\right) + f_{1}(x\left(t\right),z\left(t\right),t))\\
        & + B_{um}(\omega u\left(t\right) + f_2(x\left(t\right),z\left(t\right),t)) ,\hspace{10pt}x\left(0\right)=x_0
        \end{split}\\
        & y\left(t\right) = Cx\left(t\right)
      \end{aligned}
     \end{equation*}
 where
 \begin{equation*}
       A =
       \begin{bmatrix}
           0 & 1 & 0 & 0 & 0 & 0\\
           0 & -\frac{a_8}{a_1} & 0 & 0 & \frac{a_5}{a_1} & 0\\
           0 & 0 & 0 & 1 & 0 & 0\\
           0 & 0 & 0 & -\frac{a_{10}}{a_2} & -1.75\frac{k_c}{a_2}a_5 & \frac{a_6}{a_2}\\
           0 & 0 & 0 & 0 & -\frac{a_{16}}{a_{15}} & 0\\
           0 & 0 & 0 & 0 & 0 & -\frac{a_{18}}{a_{17}}
       \end{bmatrix},
      \end{equation*}
      
 \begin{equation*}
       B_{m} =
       \begin{bmatrix}
         0 & 0\\
         0 & 0\\
         0 & 0\\
         0 & 0\\
        \frac{a_{13}}{a_{15}} & 0\\
         0 & \frac{a_{14}}{a_{17}}
       \end{bmatrix},\\
       B_{um} =
       \begin{bmatrix}
         1 & 0 & 0 & 0\\
         0 & 1 & 0 & 0\\
         0 & 0 & 1 & 0\\
         0 & 0 & 0 & 1\\
         0 & 0 & 0 & 0\\
         0 & 0 & 0 & 0\\
       \end{bmatrix}
      \end{equation*}
      
      \begin{equation*}
       C =
       \begin{bmatrix}
        1 & 0 & 0 & 0 & 0 & 0\\
        0 & 0 & 1 & 0 & 0 & 0
       \end{bmatrix}
      \end{equation*}
      and
     \begin{equation*}
        f_{1}(x\left(t\right),t)) =
        \begin{bmatrix}
          0 & 0 & 0 & 0 & 0 & 0
        \end{bmatrix}^{\top}
     \end{equation*}
     {\scriptsize
     \begin{equation*}
        f_2(x\left(t\right),t)) =
        \begin{bmatrix}
          0 \\
          \frac{a_3}{a_1}x_5^2 -\frac{a_7}{a_1}{\rm sin}(x_1) + \frac{0.0362}{a_1}x_4^2{\rm sin}(2x_1) - \frac{a_{12}}{a_1}\big(a_3x_5^2 + a_5x_5\big)x_4{\rm cos}(x_1) \\
          0 \\
          \frac{a_4}{a_2}x_6^2 - 1.75(k_c/a_2)\big(a_3x_5^2\big) \\
          0 \\
          0
        \end{bmatrix}
     \end{equation*} }
Adaptive estimates are defined as $\hat{\theta}_1\left(t\right) \in [-50,50]\bf{1}_2$, $\hat{\theta}_2\left(t\right) \in [-50,50]\bf{1}_{4}$, $\hat{\sigma}_1\left(t\right) \in [-15,15]\bf{1}_2$, $\hat{\sigma}_2\left(t\right) \in [-15,15]\bf{1}_4$, with $\bf{1}_n:=\left[1,\ldots,1\right]^{\top}$ is a column-vector,  $\hat{\omega}_{11}\left(t\right),\hat{\omega}_{22}\left(t\right) \in [0.25,5]$, $Q = \mathbb{I}_6$, $\Gamma = 100000$ and desired poles are assigned to $-20\pm0.3i$,$-25\pm0.5i$ and $-27\pm0.5i$ and the steady state feedback gain $K$ = $10\big (\begin{smallmatrix} 1 & 0\\ 0 & 1 \end{smallmatrix}\big)$, with $k = 10$ and $D\left(s\right)=K\frac{1}{s}$. Fuzzy control parameters $k_p$, $k_d$, and $k_e$ are given as  $3.45$ , $0.05$, and $0.09$ respectively. The complete schematic of fuzzy-\Lone \,adaptive controller is shown in Fig.~\ref{Fig_5_Fuzzy_L1_5}.\\
        \begin{figure*}[ht]
           \centering
           \includegraphics[scale=0.7]{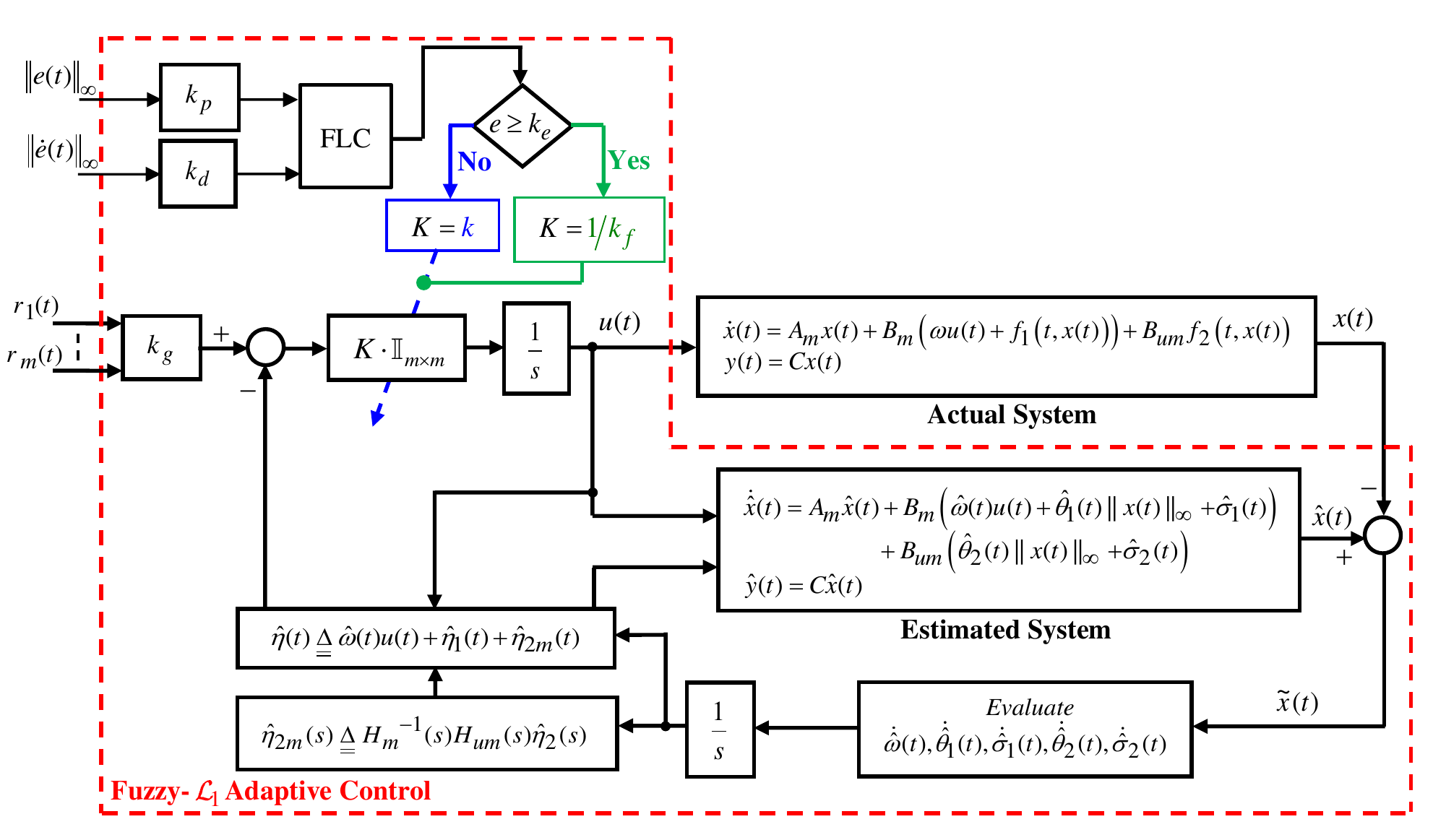}
           \caption{The proposed fuzzy-\Lone \,adaptive controller for nonlinear MIMO system.}\label{Fig_5_Fuzzy_L1_5}
        \end{figure*}
The relationship between error values and robustness margin are portrayed in Fig.~\ref{Fig_6_L1_6}, where classical \Lone \,adaptive control was applied to highly nonlinear TRMS. It can be observed that robustness margin was reduced to certain level due to the selection of \Lone \,adaptive control coefficients. Two scenarios are considered in Fig.~\ref{Fig_6_L1_6}. In the first scenario, the nonlinear TRMS was supposed to start from zero initial conditions with initial error vector $e\left(0\right) = [0,0]^{\top}$ when the reference input was $r\left(t\right) = 0.45\,{\rm sin}(0.2t)$ for both angles. The outcome of the first scenario is shown in Fig.~\ref{Fig_6_L1_6}(a) and it can be observed that excellent tracking output performance was achieved with small control signal range. In the second scenario, the initial error vector was set to $e\left(0\right) = 0.2875[1,1]^{\top}$ and desired reference was defined to be $r\left(t\right) = 0.45\,{\rm sin}(0.2t + \frac{\pi}{5})$ for both angles. In accordance with robustness margin reduction, the system stability became unstable as shown in Fig.~\ref{Fig_6_L1_6}(b). It is remarked that the observations in Fig.~\ref{Fig_6_L1_6}, \ref{Fig_6_L1_6}(a) and \ref{Fig_6_L1_6}(b) reinforces the limitations of traditional \Lone controller.\\
\indent
The aforementioned limitations are addressed by enforcing slower closed loop dynamics, which could also deteriorate the tracking performance \cite{kim_multi-criteria_2014}. These conflicting objectives are best handled using multiobjective optimization approach to achieve a compromise solution.
        \begin{figure*}[ht]
           \centering
           \includegraphics[width=15cm, height=8cm]{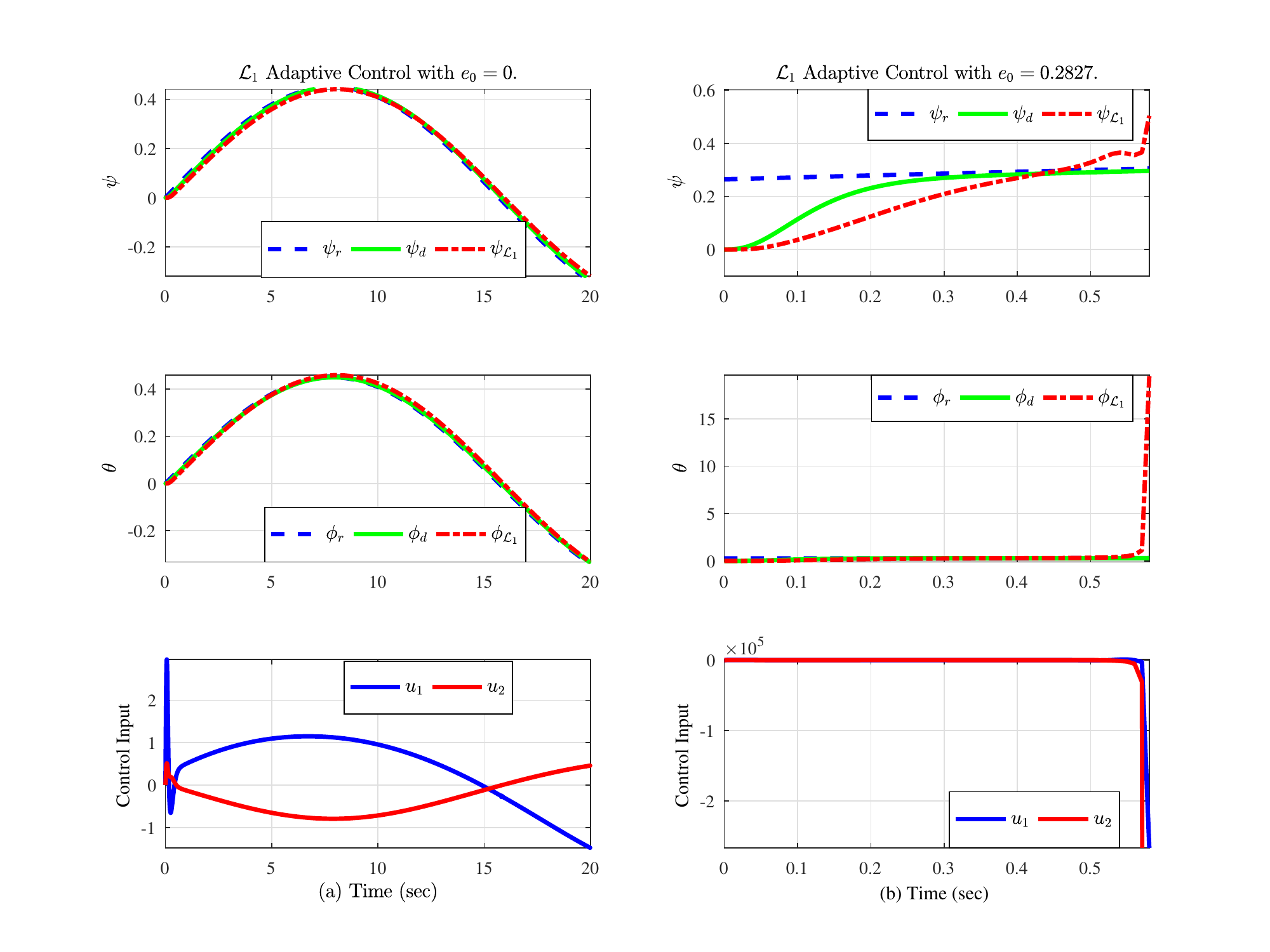}
           \caption{\Lone \,adaptive controller for nonlinear TRMS.}\label{Fig_6_L1_6}
        \end{figure*}
 \subsection{Membership Function Optimization}
Input and output membership functions for fuzzy-\Lone \,adaptive controller are here designed to improve robustness margin and reduce the control signal range in comparison with traditional \Lone \,adaptive controller. Triangular membership functions are selected and the constraint values of input and output membership functions are parameterized by lower, center and upper values. Each error-range or error-rate membership functions have four triangular linguistic variables covered by 8 optimized parameters. The output membership functions include six triangular linguistic variables covered by 16 parameters. It must be noted that the optimized parameters of membership functions are $p_i$ for $i=1,\ldots, 32$. The rule base of the proposed FLC feedback filter is defined in Table~\ref{table:Tab_Rule} with linguistic variables: $VL$ for very large, $L$ for large, $M$ for medium, $S$ for small, $VS$ for very small and $Z$ for zero.
      \begin{center}
      \begin{table}[!t]
             \setlength{\tabcolsep}{5pt}
             \setlength{\extrarowheight}{1pt}
             \caption{Rule base of FLC.} 
             \centering 
             \small
             \begin{tabular}{| c| c| c| c| c|} 
             \hline\hline 
             {\bf $\Delta e/\dot{e}$} & {\bf L} & {\bf S} & {\bf VS} & {\bf Z} \\ [0.0ex]
             \hline\hline 
             \hline
             {\bf L} & $VL$ & $VL$ & $L$ & $M$ \\[0ex]
             \hline
             {\bf S} & $VL$ & $L$ & $M$ & $S$ \\[0ex]
             \hline
             {\bf VS} & $L$ & $M$ & $S$ & $VS$ \\[0ex]
             \hline
             {\bf Z} & $M$ & $S$ & $VS$ & $Z$ \\[0ex]
             \hline\hline 
             \end{tabular}
             \label{table:Tab_Rule}
             \end{table}
      \end{center}
%
The two input membership functions have linguistic variables $L$, $S$, $VS$ and $Z$ and the output membership functions have six linguistic variables $VL$, $L$,  $M$, $S$, $VS$ and $Z$.
\subsection{MOPSO results}
Each particle in MOPSO has been designed to have 32 parameters. Parameters of $p_i, i = 1,\ldots,32$  are optimized to minimize $E$ and $U$ simultaneously. Initial settings of MOPSO algorithm are listed in Table \ref{table:Tab_PSO} with maximum number of generations equal to 50.
        \begin{table*}[ht]
        \setlength{\tabcolsep}{5pt}
        \setlength{\extrarowheight}{1pt}
        \caption{Parameters setting for MOPSO.} 
        \centering 
        \small
        \begin{tabular}{|c| c| c| c| c| c| c|} 
        \hline\hline 
        {\bf Parameter} & $\alpha$ & $c_1$ & $c_2$ & Global Set Size (Pareto-Optimal Front) & Local Set Size \\ [0.0ex]
        \hline\hline 
        {\bf Settings}  & 0.99 & 2 & 2 & 50 & 10 \\[0ex]
        \hline\hline 
        \end{tabular}
        \label{table:Tab_PSO}
        \end{table*}
The system was simulated for 23 seconds and the data was captured every 0.01 seconds. The reference input was chosen to be ${\rm cos}(0.5t)$ with zero initial conditions. After 50 generations, the optimal variables of the input and output membership functions based on the best trade-off solution are illustrated in Fig.~\ref{Fuzzy_memb_e},     \ref{Fuzzy_memb_de} and \ref{Fuzzy_memb_u}. Fig.~\ref{Fuzzy_L1_PSOObj} shows the locations of all fitness values including non-dominated solutions in MOPSO search process. MOPSO was implemented to generate a compromise solution through error and control signal range minimization. A quick observation from Fig.~\ref{Fuzzy_L1_PSOObj} reveals that reduction in control signal range increases the error with performance deterioration and vice versa. The pareto optimal front was formed by clustering a set of non-dominated solutions into best 50 solutions and most realistic compromise solution  obtained. The output performance of fuzzy-\Lone \,adaptive controller with optimized membership functions is presented in Fig.~\ref{Fuzzy_L1_out1}. The behavior of the feedback filter and the reduction in the error signal are depicted in Fig.~\ref{Fuzzy_L1_k1}.\\
        \begin{figure}[!h]
           \centering
           \includegraphics[width=7cm, height=4.5cm]{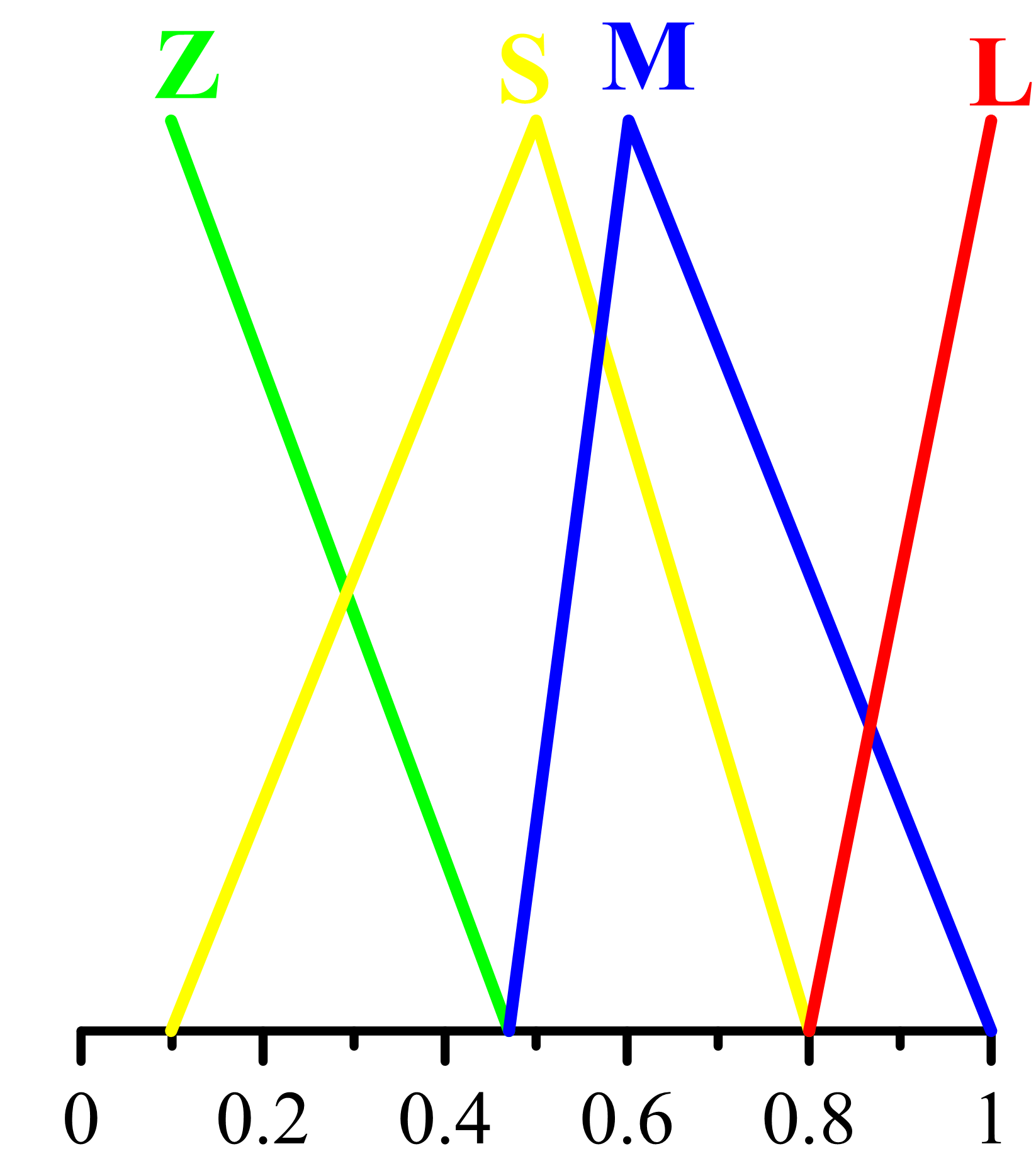}
           \caption{Optimized error membership function.}\label{Fuzzy_memb_e}
        \end{figure}
        \begin{figure}[!h]
           \centering
           \includegraphics[width=7cm, height=4.5cm]{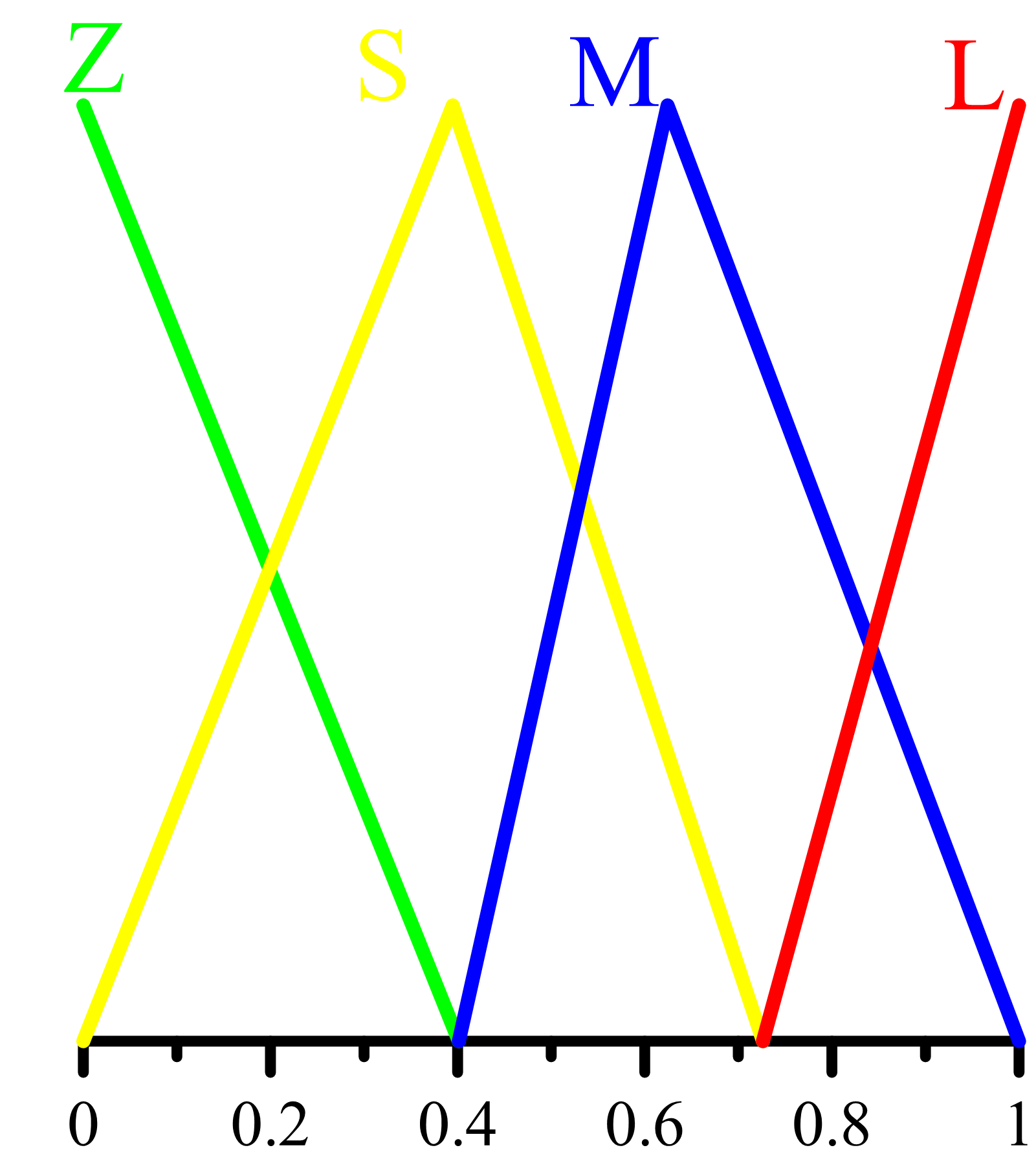}
           \caption{Optimized error rate membership function.}\label{Fuzzy_memb_de}
        \end{figure}
        \begin{figure}[!h]
           \centering
           \includegraphics[width=7cm, height=4.5cm]{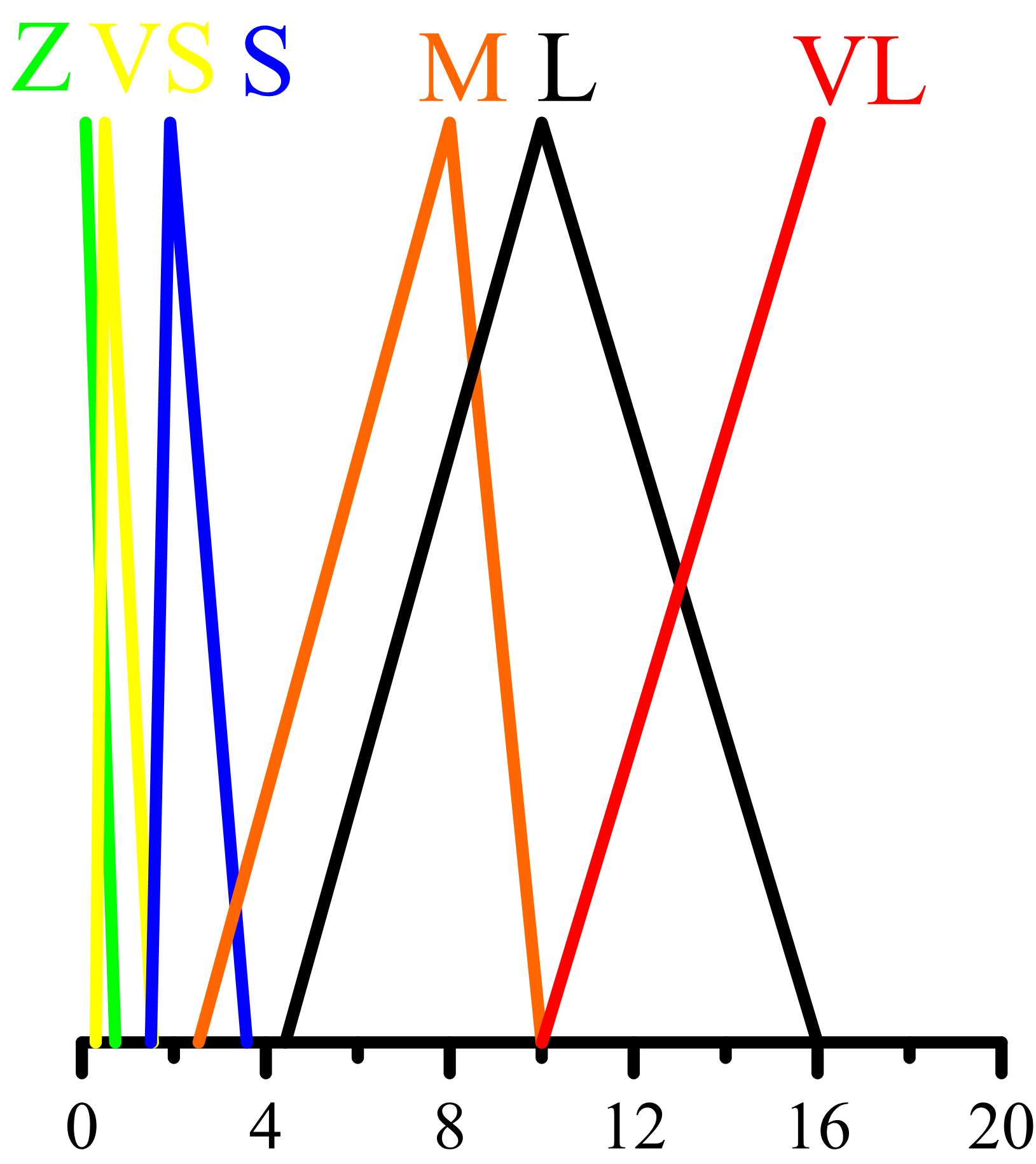}
           \caption{Optimized output membership function.}\label{Fuzzy_memb_u}
        \end{figure}

        \begin{figure*}[ht]
           \centering
           \includegraphics[width=15cm, height=8cm]{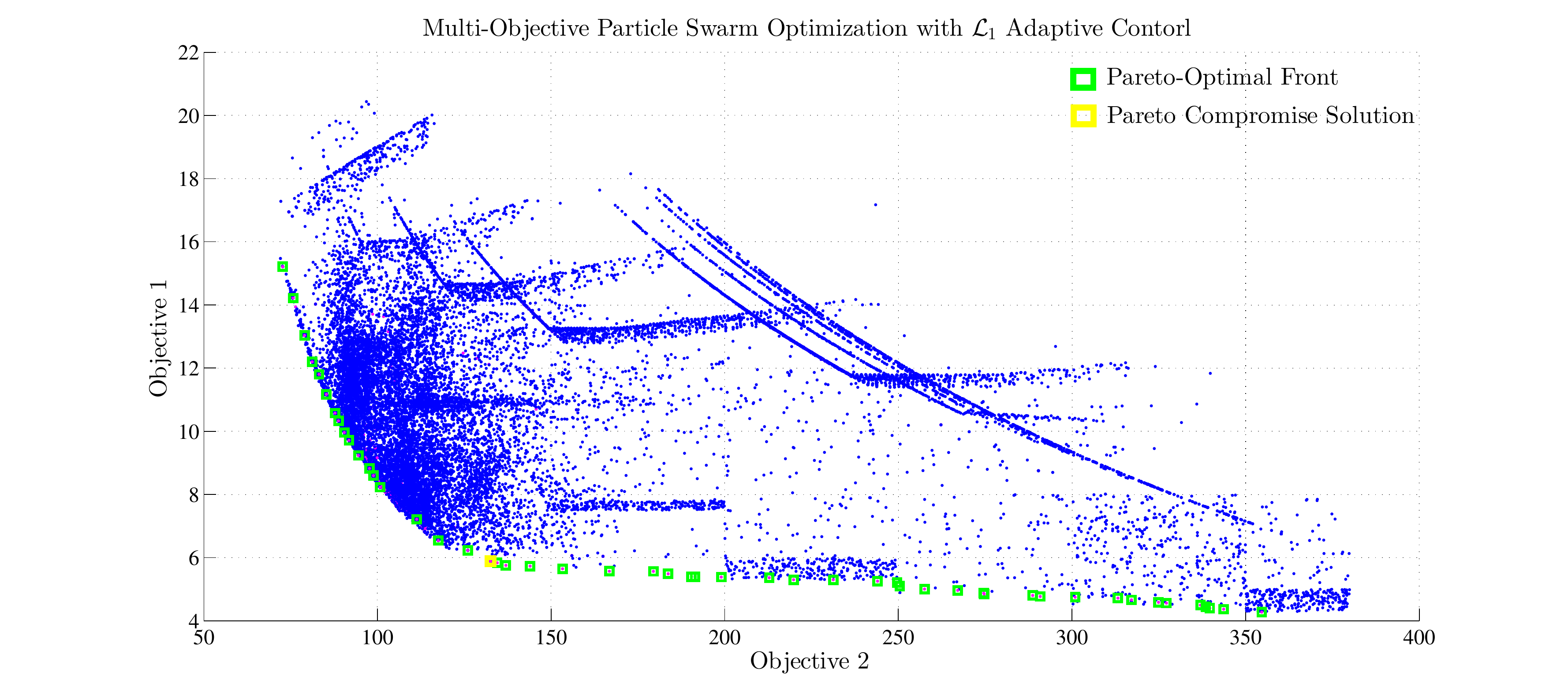}
           \caption{Multi-Objective minimization of PSO search.}\label{Fuzzy_L1_PSOObj}
        \end{figure*}
        \begin{figure*}[ht]
           \centering
           \includegraphics[scale=0.55]{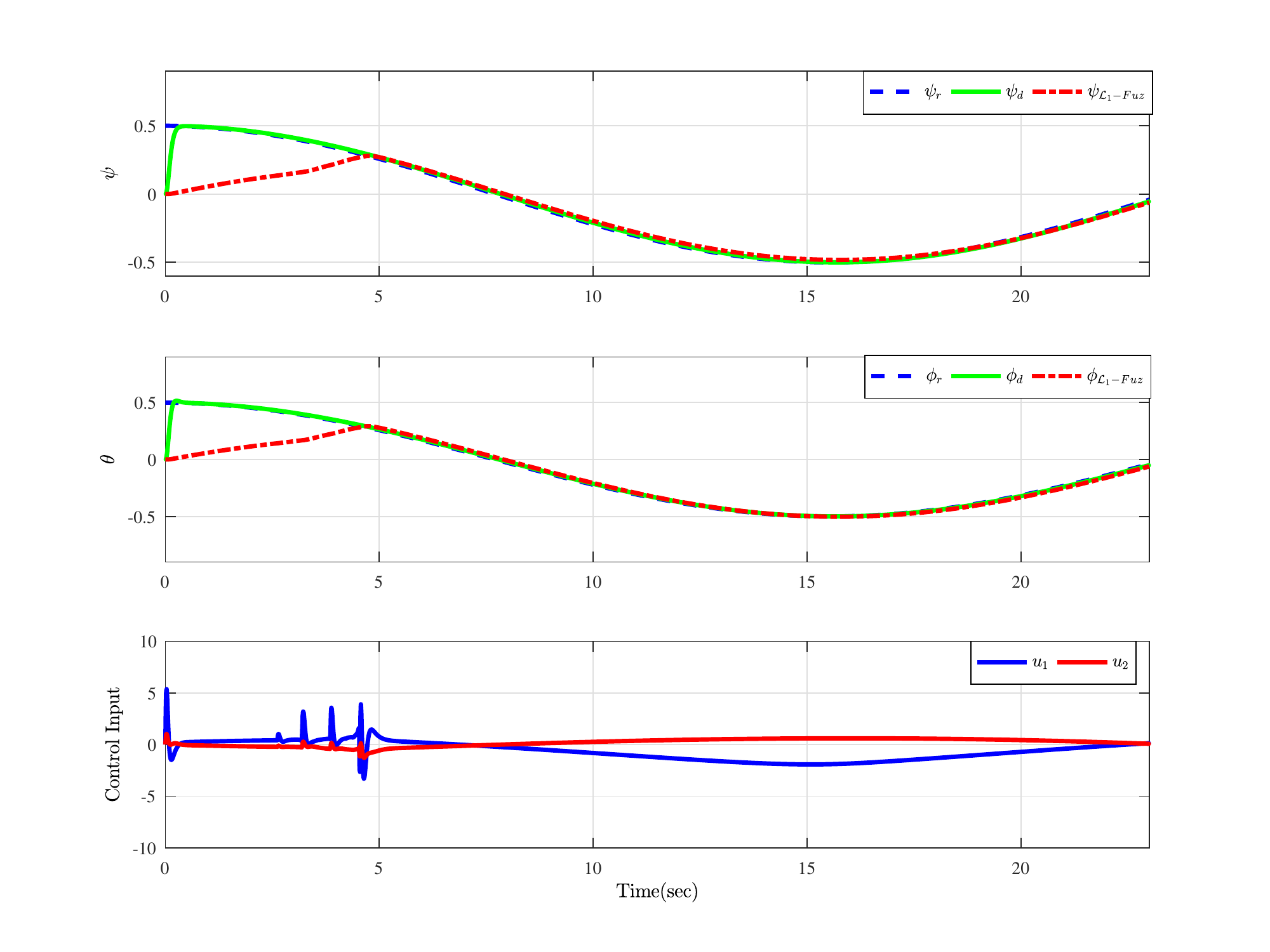}
           \caption{Output performance of fuzzy-\Lone \,adaptive controller with best compromise solution on nonlinear TRMS.}\label{Fuzzy_L1_out1}
        \end{figure*}
        \begin{figure*}[ht]
           \centering
           \includegraphics[scale=0.55]{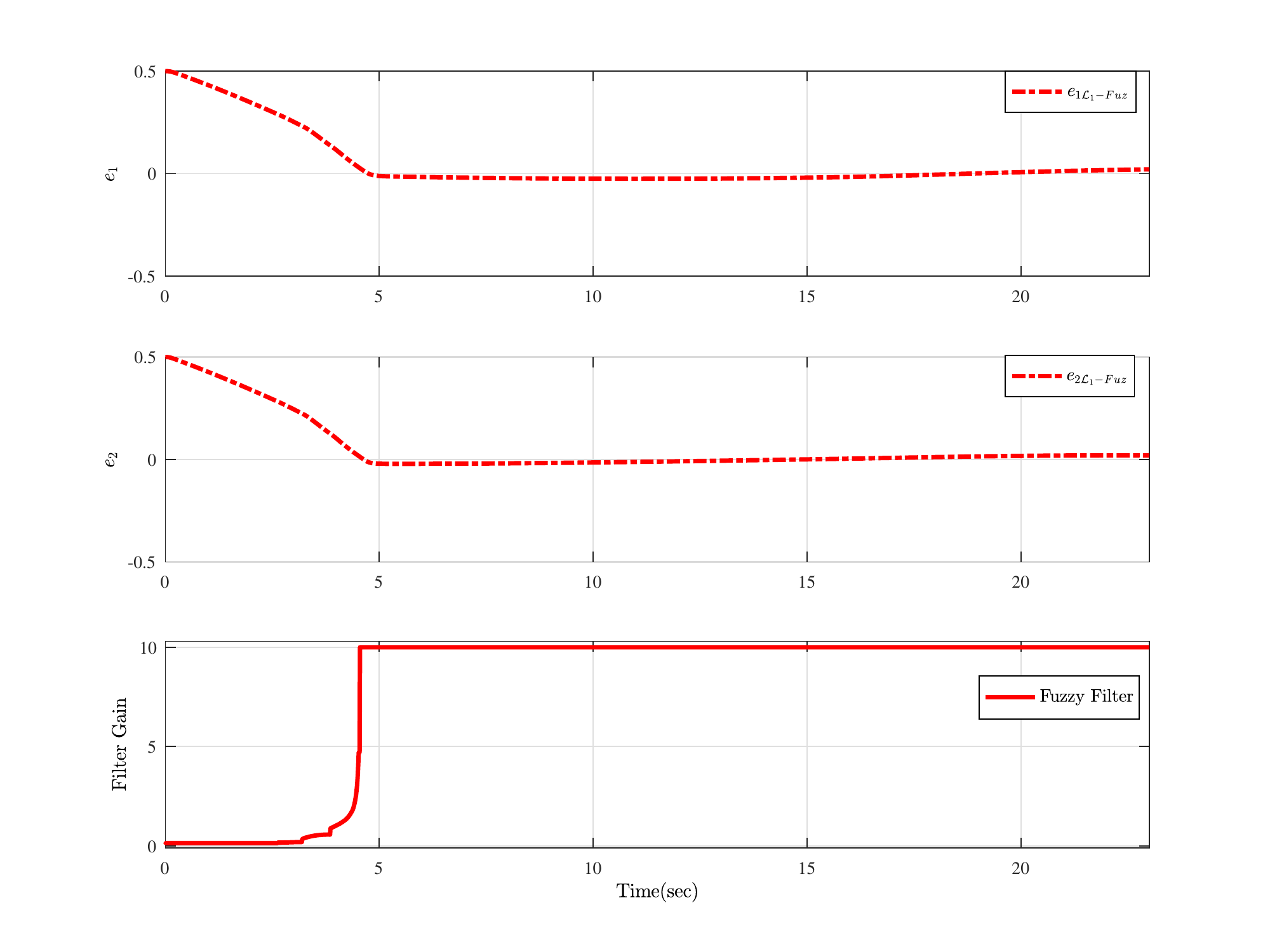}
           \caption{Performance of feedback gain filter and error for nonlinear TRMS under fuzzy-\Lone \,adaptive controller.}\label{Fuzzy_L1_k1}
        \end{figure*}
FLC with optimized parameters obtained from best compromise solution is incorporated into \Lone \,adaptive controller to tune the parameters of the feedback filter. In order to study the effect of fuzzy feedback filter with \Lone \,adaptive controller, the performance of the proposed control structure was examined for tracking capability, control signal range and robustness to uncertainties using TRMS. Two-case experiments were conducted using a composite reference signal which comprises of cosine wave, step input and other smooth functions. \\
{\bf Case 1:} the controller is implemented on a nonlinear TRMS model without uncertainties and time variant parameters. \\
{\bf Case 2:} the controller is applied to a nonlinear TRMS model with time variant and uncertain parameters.\\
For case 1, the output performance of the proposed control signal is shown in Fig.~\ref{Fuzzy_L1_out2}. It can be observed that both angles were able to track different reference signals considered. The behavior of fuzzy feedback gain filter was shown in Fig.~\ref{Fuzzy_L1_out2_2} for the purpose of simulating change in filter gain with respect to output response and error signals. Fig.~\ref{Fuzzy_L1_out2_2} presents the change in error signal with respect to tracking performance of Fig.~\ref{Fuzzy_L1_out2}.\\
        \begin{figure*}[ht]
           \centering
           \includegraphics[scale=0.55]{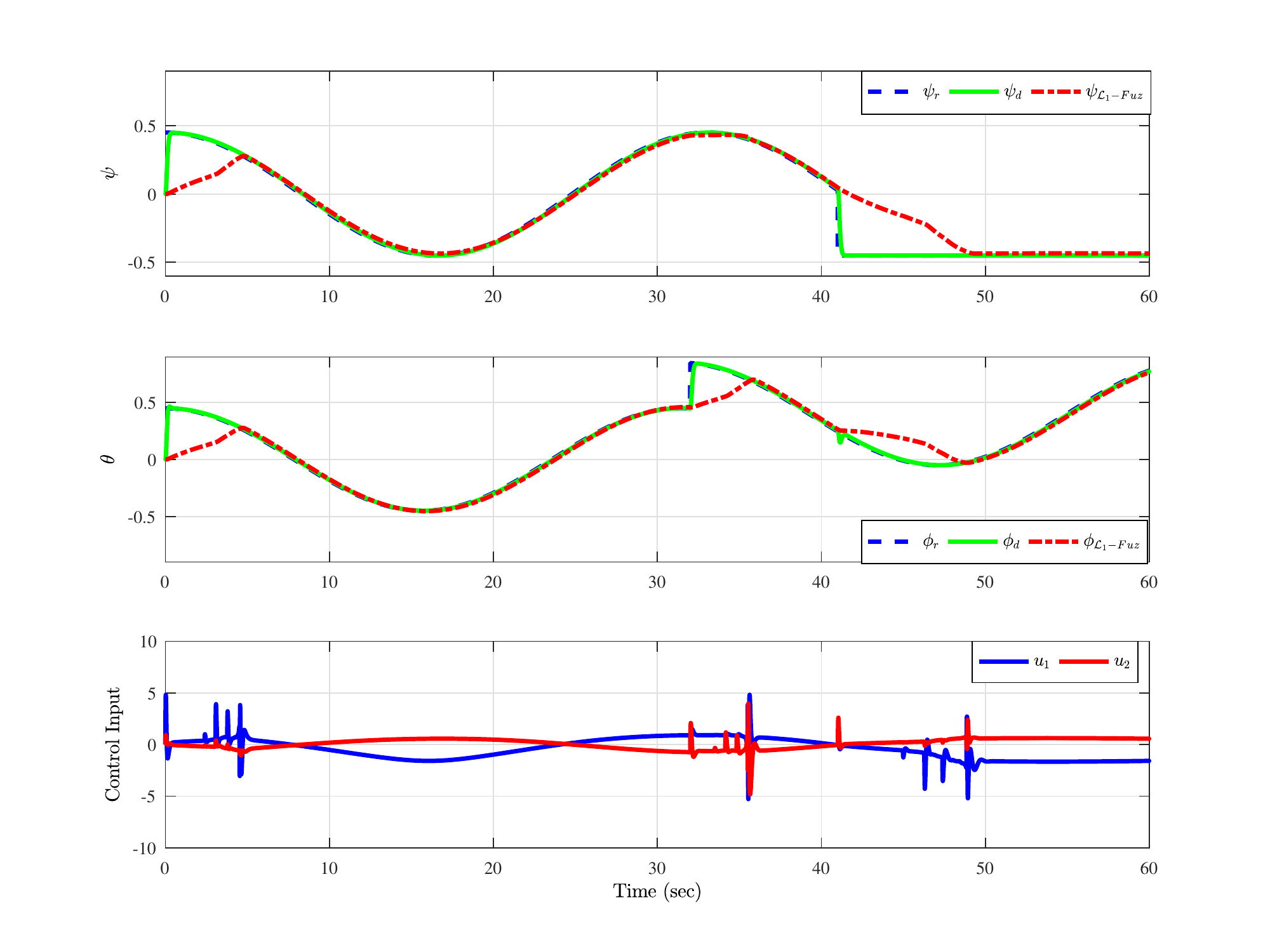}
           \caption{Case (1): Output performance of fuzzy-\Lone \,adaptive controller for nonlinear TRMS.}\label{Fuzzy_L1_out2}
        \end{figure*}

        \begin{figure*}[ht]
           \centering
           \includegraphics[scale=0.55]{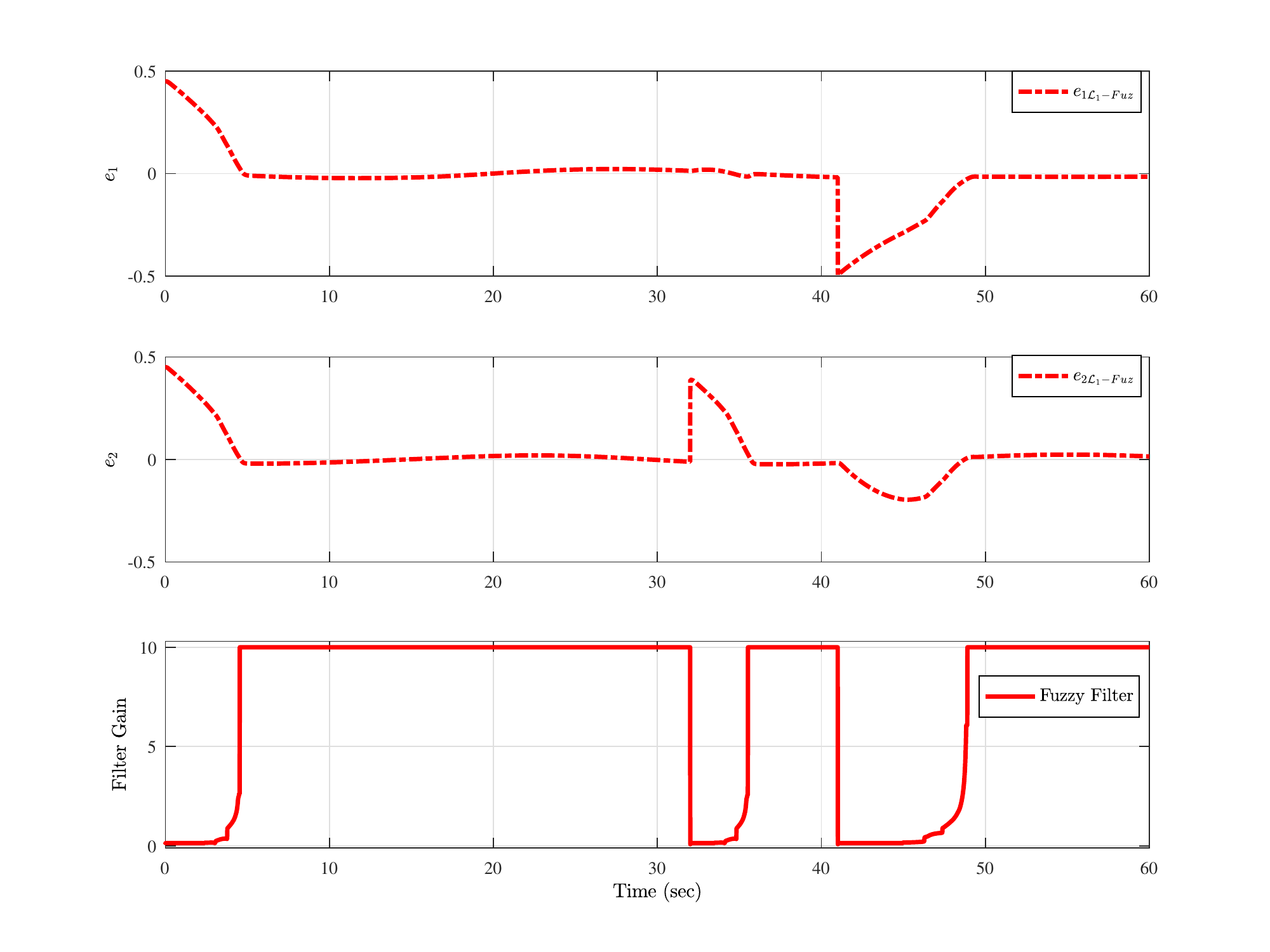}
           \caption{Case (1): Feedback gain and output error of fuzzy-\Lone \,adaptive controller for nonlinear TRMS.}\label{Fuzzy_L1_out2_2}
        \end{figure*}
 For case 2, time-varying parameters uncertainties in TRMS were considered. The output performance of fuzzy-\Lone \,adaptive controller is shown in Fig.~\ref{Fuzzy_L1_out3}. The proposed fuzzy-\Lone \,adaptive controller was applied on nonlinear TRMS with parameters $a_i = a_i (1+0.2{\rm sin}(0.3t))$ for $i=1,3,\ldots,13$ and $a_i = a_i (1+0.2{\rm cos}(0.25t))$ for $i=2,4,\ldots,14$. It is revealed in this case that the proposed approach is efficient and robust against parameter uncertainties and unmodeled dynamics. Similarly, combinations of different reference signals were considered for both pitch and yaw angles. The output performance shows that the proposed controller was able to track the desired response. The closed-loop system was also shown to be robust to presence of time-varying uncertainties as revealed in Fig.~\ref{Fuzzy_L1_out3}. Fig.~\ref{Fuzzy_L1_out3_2} shows the gain of filter response of the proposed controller and the corresponding change in errors.\\

        \begin{figure*}[ht]
           \centering
           \includegraphics[scale=0.55]{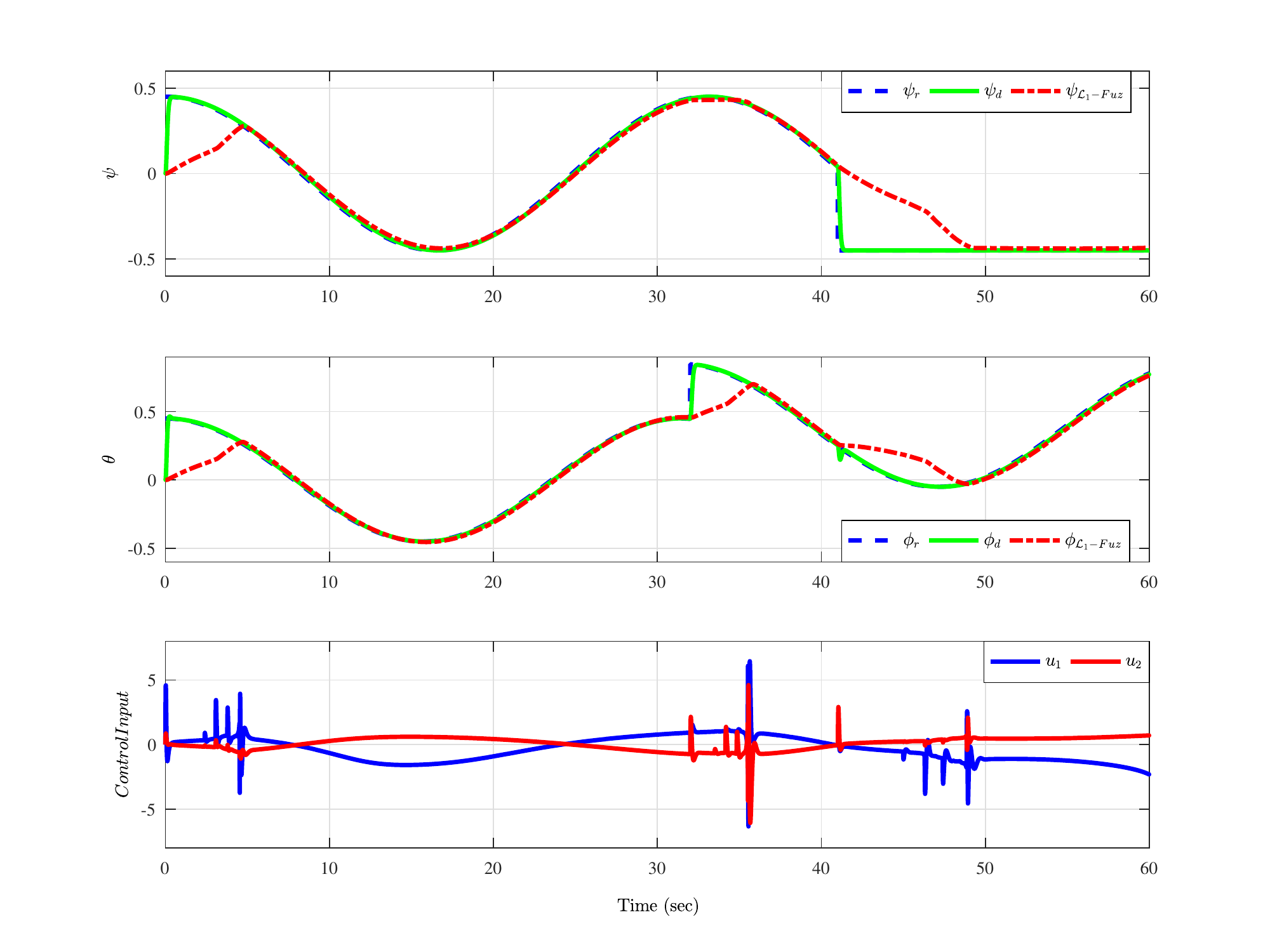}
           \caption{Case (2): Output performance of fuzzy-\Lone \,adaptive controller for nonlinear TRMS.}\label{Fuzzy_L1_out3}
        \end{figure*}

        \begin{figure*}[ht]
           \centering
           \includegraphics[scale=0.55]{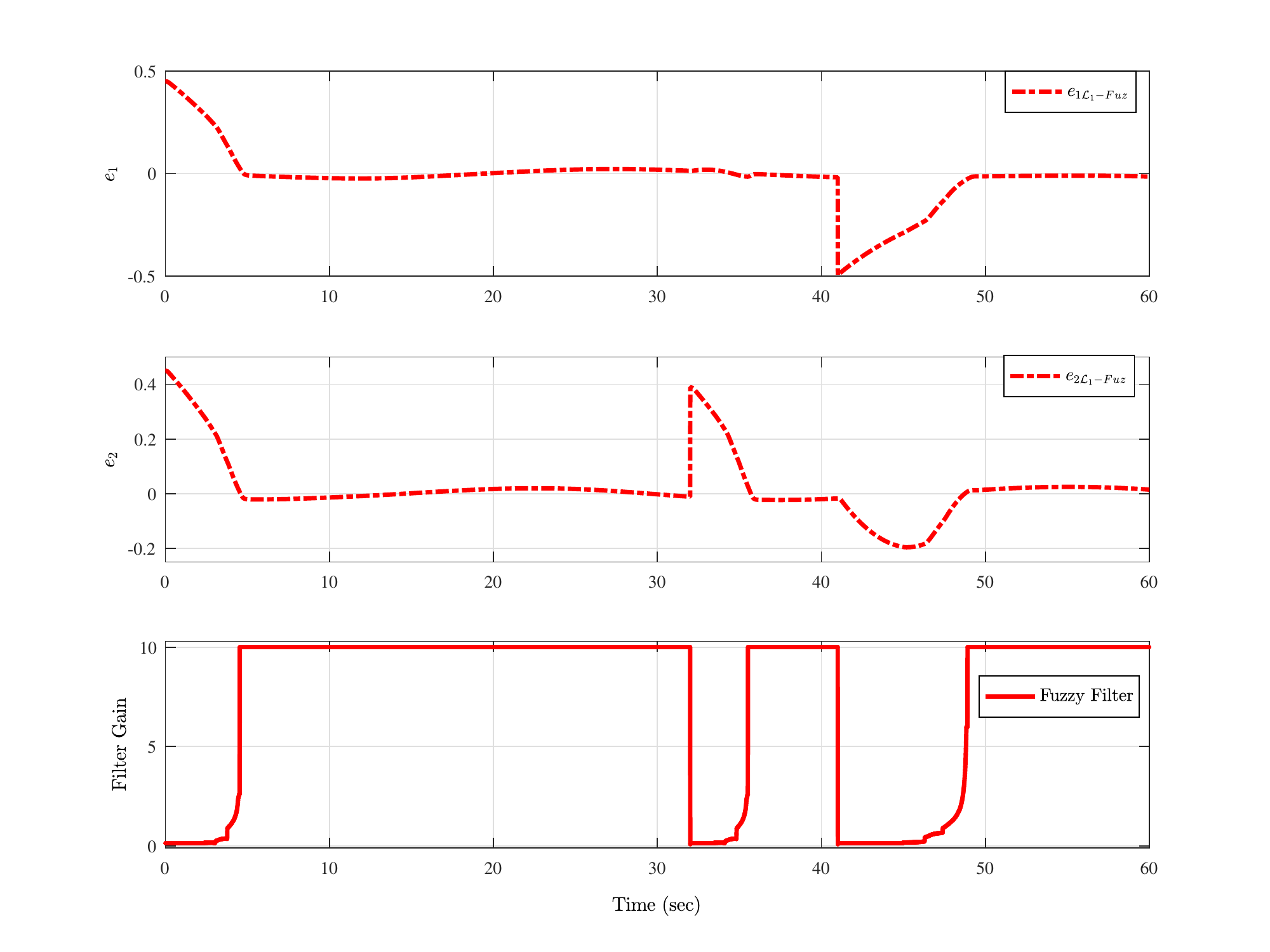}
           \caption{Case (2): Feedback gain and output error of fuzzy-\Lone \,adaptive controller for nonlinear TRMS.}\label{Fuzzy_L1_out3_2}
        \end{figure*}
It can be observed that the proposed controller guarantees smooth tracking performance and improves the robustness margins. In the case of classical \Lone \,adaptive controller, the robustness margin was exceeded and the system was into instability. The results validate the effectiveness and robustness of the fuzzy-\Lone \,adaptive controller compared to the traditional \Lone \,adaptive controller. The proposed approach is capable of tuning the feedback filter gains for both SISO and MIMO systems and providing a fast closed loop dynamics while maintaining the robustness margin and stability. The fast tracking performance and less control signal range shown are further reinforced by the results shown in Fig.~\ref{Fuzzy_L1_out2} and \ref{Fuzzy_L1_out3}.

\section{Conclusion}\label{Sec6}
In this paper, fuzzy-\Lone \,adaptive controller has been proposed for nonlinear MIMO systems. Fuzzy controller has been designed to tune the parameters of the feedback filter gain of \Lone \,adaptive controller. Multi-objective particle swarm optimization algorithm has been employed to find optimal variables for input and output membership functions based on best compromise solution between two conflicting objectives. Feedback filter parameters of the \Lone \,adaptive controller were tuned by FLC in order to improve the robustness margins. Highly nonlinear MIMO system was used to show the efficacy of the proposed approach. Results validate the effectiveness and robustness of the proposed approach on nonlinear system with time-varying uncertainties. The smooth tuning of the feedback filter enhances the robustness margin and reduces the control signal range. In addition, fast closed loop dynamics has been attained with better robustness performance.


\section*{Appendix A}\label{App_Sec6}
       \begin{table}[!h]
       \setlength{\tabcolsep}{5pt}
       \setlength{\extrarowheight}{1pt}
       \caption{Twin-Rotor Parameters.} 
       \centering 
       \small
       \begin{tabular}{|c| c| c|} 
       \hline\hline 
       {\bf variable} & {\bf Description}  & {\bf Value} \\ [0.0ex]
       \hline\hline 
       $a_1\left({\rm kg.m^2}\right)$ & Moment of inertia & $6.8\times10^{-2}$   \\ [0ex]
       $a_2\left({\rm kg.m^2}\right)$ & Moment of inertia & $2.0\times10^{-2}$  \\ [0ex]
       $a_3$ & Static parameter & $0.0135$   \\[0ex]
       $a_4$ & Static parameter & $0.0924$  \\[0ex]
       $a_5$ & Static parameter & $0.02$   \\[0ex]
       $a_6$ & Static parameter & $0.09$  \\[0ex]
       $a_7\left({\rm N.m}\right)$ & Gravity momentum & $0.32$   \\[0ex]
       $a_8\left({\rm N.m.sec/rad}\right)$ & Friction momentum & $6\times10^{-3}$   \\[0ex]
       $a_9\left({\rm N.m.sec/rad}\right)$ & Friction momentum & $1\times10^{-3}$  \\[0ex]
       $a_{10}\left({\rm N.m.sec/rad}\right)$ & Friction momentum & $0.1$   \\[0ex]
       $a_{11}\left({\rm N.m.sec/rad}\right)$ & Friction momentum & $0.01$  \\[0ex]
       $a_{12}\left({\rm rad/sec}\right)$ & gyroscopic momentum & $0.5$   \\[0ex]
       $a_{13}$ & rotor gain & $1.1$  \\[0ex]
       $a_{14}$ & rotor gain & $0.8$   \\[0ex]
       $a_{15}$ & Vertical rotor gain  & $1.1$  \\[0ex]
       $a_{16}$ & Vertical rotor gain & $1$  \\[0ex]
       $a_{17}$ & Horizontal rotor gain & $1$   \\[0ex]
       $a_{18}$ & Horizontal rotor gain & $1$  \\[0ex]

       \hline\hline 
       \end{tabular}
       \label{table:Tab_Twin}
 \end{table}
\bibliographystyle{IEEEtran}

\bibliography{Bib_L1_Fuzzy}
\end{document}